\documentclass[11pt,a4paper]{article}
\usepackage{amsmath, amssymb}
\usepackage{latexsym, color}
\usepackage{epsfig}

\numberwithin{equation}{section}

\newtheorem{thm}{Theorem}[section]

\newtheorem{lem}[thm]{Lemma}
\newtheorem{prop}[thm]{Proposition}

\def\nm{\noalign{\medskip}}

\newcommand{\qed}{\hfill \ensuremath{\square}}
\newcommand{\ds}{\displaystyle}
\newcommand{\pf}{\noindent {\sl Proof}. \ }
\newcommand{\p}{\partial}

\newcommand{\norm}[1]{\| #1 \|}

\newcommand{\eqnref}[1]{(\ref {#1})}

\newcommand{\Rbb}{\mathbb{R}}

\newcommand{\Bcal}{\mathcal{B}}

\newcommand{\Ecal}{\mathcal{E}}

\newcommand{\Ncal}{\mathcal{N}}

\newcommand{\Rcal}{\mathcal{R}}

\newcommand{\Ga}{\alpha}
\newcommand{\Gb}{\beta}
\newcommand{\Gd}{\delta}
\newcommand{\Ge}{\epsilon}

\newcommand{\Gf}{\phi}

\newcommand{\Gg}{\gamma}

\newcommand{\Gn}{\eta}

\newcommand{\Gt}{\theta}

\newcommand{\Gr}{\rho}

\newcommand{\Gs}{\sigma}

\newcommand{\Gz}{\zeta}
\newcommand{\GD}{\Delta}

\newcommand{\GG}{\Gamma}

\newcommand{\GO}{\Omega}

\newcommand{\beq}{\begin{equation}}
\newcommand{\eeq}{\end{equation}}
\newcommand{\ol}{\overline}

\begin{document}
\title{Optimal estimates of the field enhancement in presence of a bow-tie structure of perfectly conducting inclusions in two dimensions\thanks{\footnotesize This work was
supported by NRF grants No. 2015R1D1A1A01059212, 2016R1A2B4011304 and 2017R1A4A1014735.}}

\author{Hyeonbae Kang\thanks{\footnotesize Department of Mathematics, Inha University, Incheon
22212, S. Korea (hbkang@inha.ac.kr).} \and KiHyun Yun\thanks{\footnotesize Department of Mathematics, Hankuk University of Foreign Studies, Yongin-si, Gyeonggi-do 17035, S. Korea (kihyun.yun@gmail.com).}}

\maketitle

\begin{abstract}
This paper deals with the field enhancement, that is, the gradient blow-up, due to presence of a bow-tie structure of perfectly conducting inclusions in two dimensions. The bow-tie structure consists of two disjoint bounded domains which have corners with possibly different aperture angles. The domains are parts of cones near the vertices, and they are nearly touching to each other.  We characterize the field enhancement using explicit functions and, as consequences, derive optimal estimates of the gradient in terms of the distance between two inclusions and aperture angles of the corners. The estimates show that the field is enhanced beyond the corner singularities due to the interaction between two inclusions.
\end{abstract}



\section{Introduction}

This paper concerns with the field enhancement, in other words the gradient blow-up, which occurs due to presence of closely located two perfectly conducting inclusions. If we denote those inclusions by $\GO_1$ and $\GO_2$, then the problem to be considered is the following:
\beq\label{main}
\begin{cases}
\GD u = 0 \quad\mbox{in } \Rbb^2 \setminus \overline {(\GO_{1} \cup \GO_{2})}, \\
\ds u = {c}_j  \quad\mbox{on }\p \GO_{j}, \ \ j=1,2,  \\
\ds \int_{\p \GO_{1}} \p_{\nu} u ds =0, \quad j=1,2,\\
\ds u (X) - h(X) = O(|X|^{-1}) \quad\mbox{as } |X| \rightarrow \infty.
\end{cases}
\eeq
Here, $h$ is a given harmonic function in $\Rbb^2$ whose gradient represents the background field in absence of inclusions. The values $c_1$ and $c_2$ on the boundaries of inclusions are (unknown) constants depending on $h$ and inclusions, and they are determined by the zero total flux condition on $\p\GO_j$ (the third condition of \eqnref{main}).

The problem \eqnref{main} may be regarded as a conductivity problem in the context of electromagnetism, or an anti-plane elasticity problem in the context of elasticity. The solution $u$ having constant values on $\p\GO_j$ indicates that $\GO_j$ is either a perfectly conducting inclusion (the conductivity being $\infty$) or a hard inclusion (the shear modulus being $\infty$). The gradient of the solution $u$ to \eqnref{main} is either the electrical field or the stress. It can be arbitrarily large as the distance between two inclusions tends to zero, and the problem is to derive optimal estimates of the blow-up in terms of the distance between inclusions or to characterize its singular behavior using some explicit functions.

Let $u$ be the solution to \eqnref{main}, and let
\beq
\Ge:= \mbox{dist} (\GO_1, \GO_2).
\eeq
The problem of estimating $|\nabla u|$ in terms of $\Ge$ was first raised in \cite{bab} in relation to stress analysis of composites, and there has been significant progress on this problem in the last decade or so. It is proved that the optimal blow-up rate of $|\nabla u|$ is $\Ge^{-1/2}$ in two dimensions \cite{AKL-MA-05, Y}, and $(\Ge|\ln\Ge|)^{-1}$ in three dimensions \cite{BLY-ARMA-09}. The singular behavior of $\nabla u$ is characterized asymptotically in \cite{ACKLY-ARMA-13, KLY-MA-15, KLY-JMPA-13, KLY-SIAP-14}.
It is worth mentioning that the gradient estimate was extended to the case of insulating inclusions \cite{AKL-MA-05, BLY-CPDE-10, Yun-JDE-16}, non-homogeneous case \cite{DL-arXiv}, $p$-Laplacian \cite{GN-MMS-12}, and the Lam\'e system of the linear elasticity \cite{BLL-ARMA-15, KY-arXiv}.

\begin{figure}[h!]
\begin{center}
\epsfig{figure=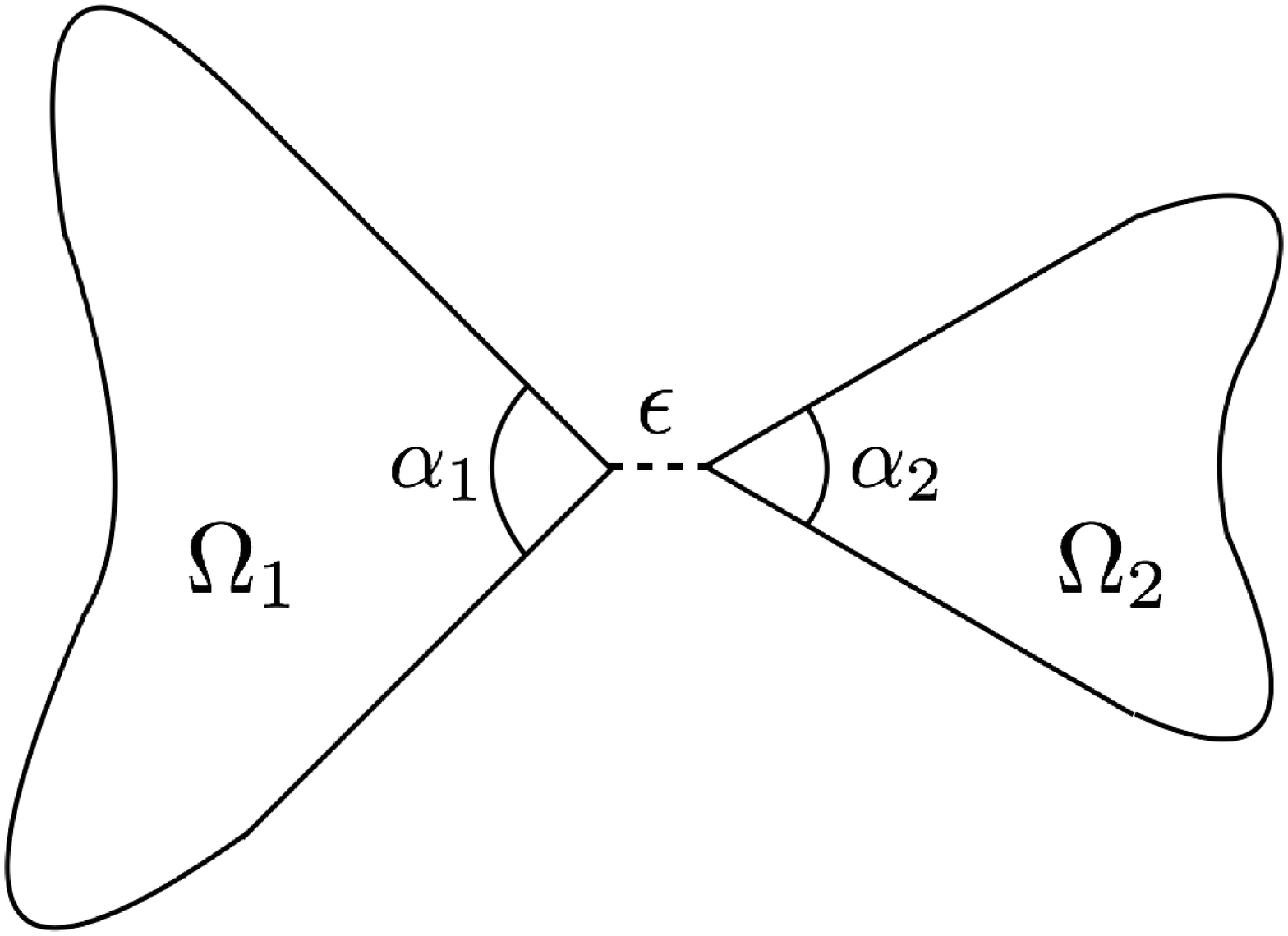,width=5cm}
\end{center}
\caption{A bow-tie structure}\label{Fig1}
\end{figure}

All the above mentioned work deal with inclusions with smooth boundaries, $C^{2,\Ga}$ boundaries to be precise. However, for some purposes such as imaging, the stronger gradient blow-up is desirable (see, for example, \cite{PBFLN}), and for that it is natural to consider structures consisting of inclusions with corners instead of those with smooth boundaries. One such a structure is the bow-tie structure as depicted in Figure \ref{Fig1}. In fact, there are some works dealing with field enhancement on the bow-tie structure, for which we refer to the above mentioned article and references therein. However, a rigorous quantitative analysis of the field enhancement is missing as far as we are aware of, and it is the purpose of this paper to derive optimal estimates for the field enhancement due to presence of the bow-tie structure.

The gradient blow-up due to presence of the bow-tie structure is particularly interesting since two different types of singularities of the gradient are expected: one due to presence of corners and the other due to the interaction between closely located inclusions. Since the fundamental work of Kontratiev \cite{Kondra-TMMS-67} on the corner singularity (see also \cite{Grisvard-book, KMR-book}), it is well known that near the vertex of a corner with the aperture angle $\Ga$, the gradient of the solution to the Laplace equation has singularity of the form $|X-V|^{-1+\Gb}$ where $\Gb = \frac{\pi}{2\pi - \Ga}$. In the bow-tie structure, two vertices are closely located. So one can expect the corner singularity may be amplified due to the interaction of the inclusions. The results of this paper confirm the expectation in a quantitatively precise manner.

Let $V_j$ be the vertex of $\GO_j$ and the aperture angle at $V_j$ be $\Ga_j$ ($j=1,2$). We show, for example, that the gradient of the solution to \eqnref{main} behaves near the vertex $V_j$ as
\beq\label{blow-up1}
\frac{\Ge^{-\Gb_j}}{|\log \Ge|} |X-V_j|^{-1+\Gb_j},
\eeq
where
\beq\label{Gbdef}
\Gb_j := \frac{\pi}{2\pi - \Ga_j}, \quad j=1,2,
\eeq
(see Theorem \ref{1st_cor}). So, the corner singularity $|X-V_j|^{-1+\Gb_j}$ is amplified by the factor of $\Ge^{-\Gb_j}/|\log \Ge|$. According to \eqnref{blow-up1}, the magnitude of the gradient is larger than $\Ge^{-1}/|\log \Ge|$ provided that $|X-V_j| \le \Ge$. This blow-up rate is much larger than the one in the case of inclusions with smooth boundaries, where the blow-up rate is $\Ge^{-1/2}$ as mentioned before. We also obtain estimates of $|\nabla u(X)|$ when $X$ is relatively away from the vertices (Theorem \ref{cortwo}), which shows as well that the corner singularity is amplified on the bow-tie structure.

In this paper we assume that the aperture angles $\Ga_1$ and $\Ga_2$ satisfies $0< \Ga_j \le \pi$ and
\beq\label{anglecond}
C_1 \le \Ga_1+\Ga_2 \le 2\pi - C_2
\eeq
for some positive constants $C_1$ and $C_2$. Estimates to be derived in the sequel depend on $C_1$ and $C_2$. So, one of angles is allowed to be $\pi$.

This paper is organized as follows. In section \ref{sec:corner} we define the bow-tie structure and investigate the corner singularity in the case of the present paper, where there are two corners. In section \ref{sec:cone} we introduce some auxiliary functions to be used in later sections. In section \ref{sec:singular}, we introduce singular functions which capture singular behavior of $\nabla u$. The main results and their proofs are presented in section \ref{sec:main}.  The last section is to estimate a coefficient, which is crucial in achieving lower bounds.

Throughout this paper, the expression $A \lesssim B$ implies that there is a constant $C$ independent of $\Ge$ and the background potential $h$ such that $A \le CB$, and $A \simeq B$ implies that both $A \lesssim B$ and $B \lesssim A$ hold.

\section{Corner singularities}\label{sec:corner}

Let $\GG_1$ be the open cone in the left-half space with the vertex at $S_1:=(-1/2,0)$ and the aperture angle $\Ga_1$ at $S_1$, namely,
\beq
\GG_1= \{ X=(x,y) ~:~ |y| < -\tan(\Ga_1/2) (x+1/2) \},
\eeq
and let
\beq
\GG_2= \{ X=(x,y) ~:~ |y| < \tan(\Ga_2/2) (x-1/2) \},
\eeq
which is the open cone in the right-half space with the vertex at $S_2:=(1/2,0)$ and the aperture angle $\Ga_2$ at $S_2$. Then inclusions $\GO_1$ and $\GO_2$, which constitute the bow-tie structure, can be defined locally as translates of $\GG_1$ and $\GG_2$: Let
\beq\label{LGe}
L_\Ge:= \frac{1}{2}(-\Ge +1) \quad\mbox{and}\quad R_\Ge:= \frac{1}{2}(\Ge-1).
\eeq
Then, $\GO_1$ and $\GO_2$ are bounded domains whose boundaries are smooth except at vertices $V_1:=(-\Ge/2,0)$ and $V_2:=(\Ge/2,0)$, respectively, such that
\beq\label{BGd}
\GO_1 \cap B_\Gd= (\GG_1+L_\Ge) \cap B_\Gd \quad\mbox{and}\quad \GO_2 \cap B_\Gd= (\GG_2+R_\Ge) \cap B_\Gd
\eeq
for some $\Gd>0$. Here and throughout this paper $B_\Gd(P)$ denotes the open disk of radius $\Gd$ centered at $P$, and if $P=O$, the origin, then we simply denote it by $B_\Gd$. We emphasize that if $\GG_1$ and $\GG_2$ are scaled by $\Ge$, then they become $\GO_1$ and $\GO_2$ near the origin $O$, that is,
$$
\Ge (\GG_1 \cup \GG_2) \cap B_\Gd = (\GO_1 \cup \GO_2) \cap B_\Gd.
$$

\begin{figure}[h!]
\begin{center}
\epsfig{figure=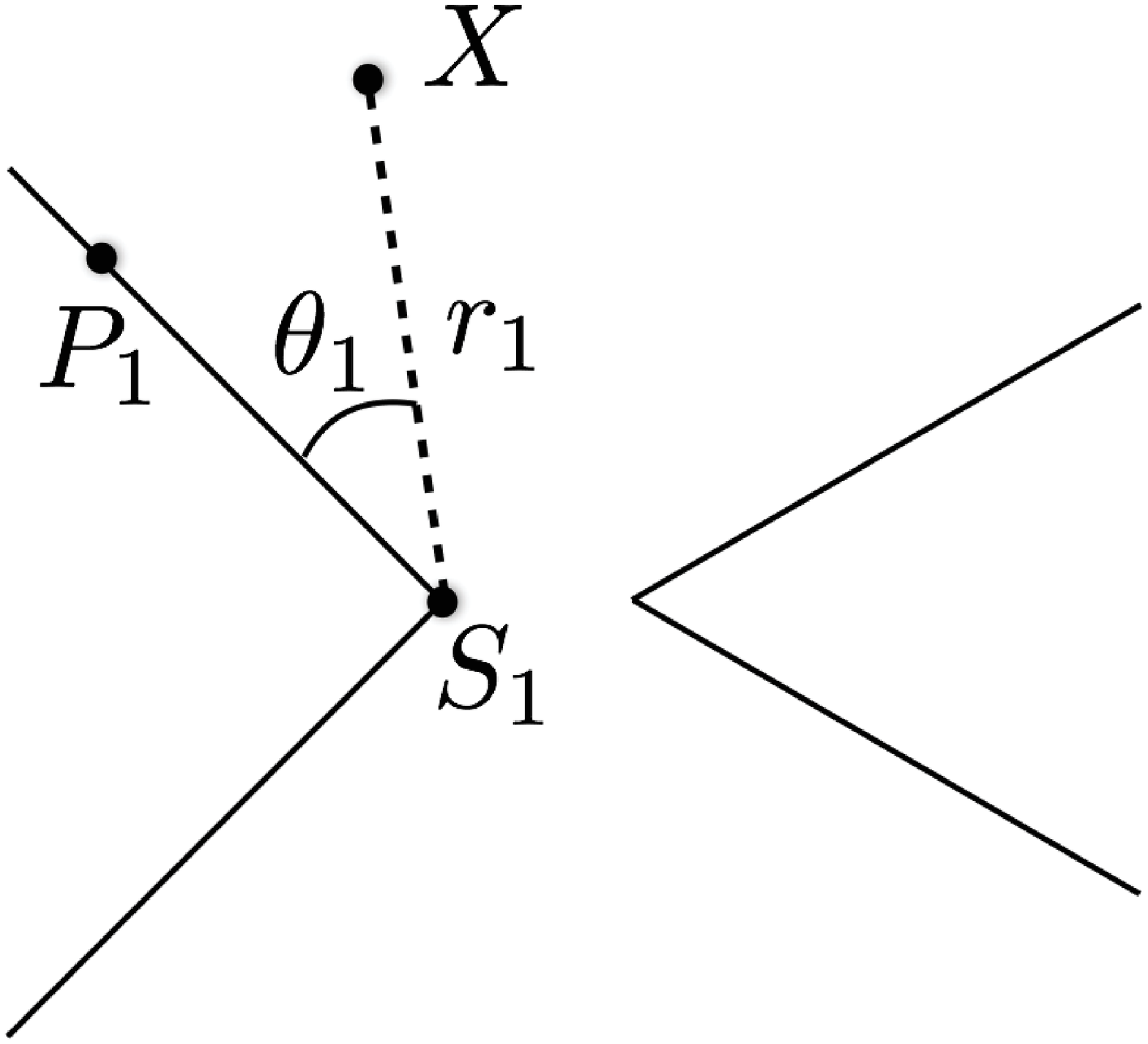,width=6cm}
\end{center}
\caption{Polar coordinates with respect to vertices}\label{Fig2}
\end{figure}

Choose and fix a point $P_j$ on $\p\GG_j \cap \{ (x,y): y>0 \}$ for $j=1,2$, and define $\Gt_j(X)$ and $r_j(X)$ for $X \in \Rbb^2 \setminus \overline {\GG_j}$ ($j=1,2$) by
\beq\label{angle}
\mbox{$\Gt_j(X)=$the angle between $\overrightarrow {S_j P_j}$ and $\overrightarrow{S_j X}$},
\eeq
and
\beq\label{rldef}
r_j(X) =  |X-S_j|.
\eeq
See Figure \ref{Fig2}. Let $\Gb_j$ be the number defined by \eqnref{Gbdef}
and define the functions $\Bcal_j(X)$ for $j=1,2$ by
\begin{align}
\Bcal_j(X) &= r_j (X)^{\Gb_j} \sin (\Gb_j \Gt_j (X)), \quad X \in \Rbb^2 \setminus \overline {\GG_j} \,  .
\end{align}
The function $\Bcal_j$ is the singular part of the solution to the elliptic problem when the domain has a corner of the angle $2\pi-\Ga_j$ (see \cite{Grisvard-book}). Note that $\nabla \Bcal_j$ has a singularity of order $\Gb_j-1$ at the vertex $S_j$. In fact, one can easily see that
\beq\label{naBcalj}
|\nabla \Bcal_j(X)| = \Gb_j |X-S_j|^{\Gb_j-1}, \quad j=1,2.
\eeq

Since the bow-tie structure has two corners, the function $\Bcal_1 - b\Bcal_2$ plays an important role in the sequel, where $b$ is a positive number determined by angles $\Ga_1$ and $\Ga_2$ (see Theorem \ref{1st_thm}). In the following we look into its properties.

\begin{lem}\label{One}
For $j=1,2$, let $\nu_j$ be the normal vector of $\p \GG_j$ pointing towards the inside of $\GG_j$. The following holds:
\begin{itemize}
\item[(i)] for $X \in \p\GG_1 \setminus \{ S_1 \}$
\beq\label{Gammaone}
\p_{\nu_1} \Bcal_1(X) = - | \nabla \Bcal_1(X)| < 0 \quad\mbox{and}\quad  \p_{\nu_1} \Bcal_2\left ( X\right) \ge 0,
\eeq

\item[(ii)] for $X \in \p\GG_2 \setminus \{ S_2 \}$
\beq
\p_{\nu_2} \Bcal_2\left ( X\right) = - \left| \nabla \Bcal_2\left ( X\right)\right | < 0 \quad\mbox{and}\quad  \p_{\nu_2} \Bcal_1\left ( X\right) \ge 0.
\eeq
\end{itemize}
\end{lem}
\pf
Let us prove (i) only. (ii) can be proved similarly.

Since $\Bcal_1=0$ on $\p\GG_1$, we have $\p_{\nu_1} \Bcal_1(X) \nu_1 = \nabla \Bcal_1(X)$ for $X \in \p\GG_1 \setminus \{ S_1 \}$. Moreover, since $\Bcal_1 >0 $ in $\Rbb^2 \setminus \overline{\GG_1}$, $\p_{\nu_1} \Bcal_1(X) <0$ by Hopf's lemma. So, the first inequality in \eqnref{Gammaone} follows.

To prove the second inequality in \eqnref{Gammaone}, we note that $X=(x, y)$ in $\Rbb^2 \setminus \GG_2$ can be represented, using the polar coordinates with respect to $S_2$ as defined in \eqnref{angle} and \eqnref{rldef}, as
$$
x-\frac 1 2 =r_2 \cos \left( \Gt_2 + \frac{\Ga_2}{2} \right), \quad y=r_2 \sin \left( \Gt_2 + \frac{\Ga_2}{2} \right).
$$
So we see that
$$
U_{-\Gt-\Ga_2/2} \nabla \Bcal_2(X)= (\p_{r_2} \Bcal_2(X), r_2^{-1} \p_{\Gt_2} \Bcal_2(X))^T \ \ (T \mbox{ for transpose}),
$$
where $U_{-\Gt_2-\Ga_2/2}$ is the matrix of rotation by $-\Gt_2-\Ga_2/2$ with respect to $S_2$. So, one can see that
\beq\label{nabcal2}
\nabla \Bcal_2(X) = \Gb_2 r_2^{\Gb_2-1} \left( -\sin \left( \Gt_2+ \frac{\Ga_2}{2} - \Gb_2 \Gt_2 \right), \cos \left( \Gt_2+ \frac{\Ga_2}{2} - \Gb_2 \Gt_2 \right) \right)^T.
\eeq
One can also see that
$$
\nu_1= -\left( \sin \left( \frac{\Ga_1}{2} \right), \cos \left( \frac{\Ga_1}{2} \right) \right)^T.
$$
It then follows that
\beq\label{Gnonecdot}
\nu_1 \cdot \nabla \Bcal_2(X)= - \cos \left( (1-\Gb_2)\Gt_2 + \frac{\Ga_1+\Ga_2}{2} \right).
\eeq

If $X =(x,y) \in \p\GG_1$ with $y \ge 0$, then we have
$$
\pi - \frac{\Ga_1+\Ga_2}{2} < \Gt_2(X) \le \pi - \frac{\Ga_2}{2},
$$
and hence
\begin{align*}
(1-\Gb_2) \left( \pi - \frac{\Ga_1+\Ga_2}{2} \right) + \frac{\Ga_1+\Ga_2}{2} &< (1-\Gb_2)\Gt_2 + \frac{\Ga_1+\Ga_2}{2} \\
&\le (1-\Gb_2) \left( \pi - \frac{\Ga_2}{2} \right) + \frac{\Ga_1+\Ga_2}{2}.
\end{align*}
Note that the term on the far left-hand side is linear in $\Ga_1$. By substituting $\Ga_1=0$ and $\pi$, we see that it is bounded from below by $\pi/2$. Likewise one can see that the term on the far right-hand side is bounded from above by $\pi$. So, we have
$$
\frac{\pi}{2} < (1-\Gb_2)\Gt_2 + \frac{\Ga_1+\Ga_2}{2} < \pi.
$$
Now the second inequality in \eqnref{Gammaone} follows from \eqnref{Gnonecdot}, and the proof is complete for $X=(x,y)$ with $y > 0$. The case when $y<0$ can be treated in the same way. This completes the proof. \qed

The following proposition shows that there is no cancellation between two gradients $\nabla\Bcal_1$ and $-\nabla\Bcal_2$ entering $\nabla \Bcal_1 - b \nabla \Bcal_2$.

\begin{prop}\label{A} Suppose that $b>0$. The following holds:
\beq\label{B1B2est}
\left|\nabla \Bcal_1 (X) \right| + \left| \nabla \Bcal_2 (X) \right| \lesssim \left|\nabla \Bcal_1 (X) - b \nabla \Bcal_2 (X) \right|
\eeq
for  $X \in \Rbb^2  \setminus \overline{\GG_1 \cup \GG_2}$.
\end{prop}
\pf
We begin the proof by showing that \eqnref{B1B2est} holds on $\p \left( B_R \setminus (\GG_1 \cup \GG_2) \right) \setminus \{ S_1, S_2 \}$ for any sufficiently large $R >0$. To do so, we define $l_R^\pm$ by
$$
l_R^{\pm} := \left \{ X ~:~  X=(x,y) \in \p B_R \setminus (\GG _1 \cup \GG_2 ), \ \pm y > 0  \right \},
$$
respectively, so that we have
\begin{align*}
&\p \left( B_R \setminus (\GG _1 \cup \GG_2 ) \right) \setminus \left \{ S_1, S_2 \right\} \\
&= l_R^{+} \cup l_R^{-} \cup ( ( \p \GG_1 \cap B_R ) \setminus \{S_1 \}) \cup ( ( \p \GG_2 \cap B_R ) \setminus \{S_2 \}).
\end{align*}

First, we prove \eqnref{B1B2est} on $l_R^+$.  Let $(r, \Gt)$ be the usual polar coordinates with respect to $O$, the origin. We write $X = r e^{i\Gt}$ by identifying points in $\Rbb^2$ with complex numbers. Note that
\beq\label{loglog}
\log \left(R e^{i\Gt} - \frac 1 2 \right) = \log r_2 + i \left(\theta_2 + \frac {\Ga_2} 2 \right)
\eeq
on $l_R^+$. Since $|R e^{i\Gt} - 1/2| = r_2$, we infer that
\beq\label{r2/R}
\frac{r_2}{R} = 1 + O\left(\frac 1 R\right)
\eeq
for all sufficiently large $R$. By differentiating both sides of \eqnref{loglog} with respect to $\Gt$, we obtain
$$
\frac{Rie^{i\Gt}}{Re^{i\Gt}-\frac{1}{2}} = \frac{1}{r_2} \frac{\p r_2}{\p\Gt} + i \frac{\p \Gt_2}{\p\Gt},
$$
from which we infer that
$$
\left| \frac {\p \Gt_2} {\p \Gt}  - 1  \right| + \frac 1 {r_2}\left| \frac {\p r_2} {\p \Gt} \right| \leq \frac C R \quad \mbox{on } l_R^+,
$$
provided that $R$ is sufficiently large. Here and in the following the constant $C$ may differ at each occurrence. It then follows that
\begin{align}
\frac {\p \Bcal_2 (X)} {\p \Gt} &= \frac {\p \Gt_2} {\p \Gt} \frac {\p \Bcal_2 (X) } {\p \Gt_2}  + \frac {\p r_2} {\p \Gt} \frac {\p \Bcal_2 (X)} {\p r_2} \notag\\&= \left(1+ O\left(\frac 1 R\right)\right) \frac {\p \Bcal_2 (X) } {\p \Gt_2} + O\left(\frac 1 R\right) r_2 \frac {\p \Bcal_2 (X)} {\p r_2} \label{deri_B_2}
\end{align}
for all sufficiently large $R$.

One can see that there is $R_0$ such that $0 \le \Gb_j \Gt_2(X) \le \pi/2 - \Gt_0$ for some $\Gt_0 >0$ for all $X \in l_R^+$ with $R \ge R_0$. So we have
$$
\sin (\Gb_2 \Gt_2 (X)) \le C \cos (\Gb_2 \Gt_2 (X)),
$$
and hence
\beq\label{GtGt2}
0< r_2 \frac{\p \Bcal_2 (X)}{\p r_2} \leq C \frac {\p \Bcal_2 (X)} {\p \Gt_2} .
\eeq
This together with \eqnref{deri_B_2} shows that there is $R_1$ such that
\beq\label{deri_B_2two}
\frac {\p \Bcal_2(X)}{\p \Gt} \simeq \frac {\p \Bcal_2 (X) }{\p \Gt_2}
\eeq
for all $X \in l_R^+$ with $R \ge R_1$. We also have
\beq\label{deri_B_2three}
|\nabla \Bcal_2 (X) | \leq C \left( \frac 1{r_2}  \frac{\p \Bcal_2(X)}{\p \Gt_2} + \frac{\p \Bcal_2 (X)}{\p r_2} \right) \leq C \frac 1{r_2}  \frac{\p \Bcal_2(X)}{\p \Gt_2}.
\eeq
It then follows from \eqnref{r2/R}, \eqnref{deri_B_2two} and \eqnref{deri_B_2three} that
\beq\label{deri_B_2four}
|\nabla \Bcal_2 (X) | \leq C \frac 1{R}  \frac{\p \Bcal_2 (X)}{\p \Gt}, \quad X\in l_R ^+ ,
\eeq
for all $R\ge R_1$.

In the same way, one can show that there is a positive constant $C$ independent of $R$ such that
$$
\frac 1 {r_1} \frac{\p \Bcal_1(X)}{\p\Gt_1} \leq -C \frac 1 R \frac{\p \Bcal_1(X)}{\p \Gt}
$$
and
\beq\label{deri_B_2five}
|\nabla \Bcal_1 (X) |  \leq - C \frac 1 R \frac{\p \Bcal_1(X)}{\p \Gt}, \quad X \in l_R ^+ ,
\eeq
provided that $R \ge R_1$.

Combining \eqnref{deri_B_2four} and \eqnref{deri_B_2five}, we arrive at
\begin{align}
|\nabla \Bcal_1 (X)| + |\nabla \Bcal_2 (X)| & \leq - C \left( \frac 1 R \frac {\p \Bcal_1 (X)} {\p \Gt} - b \frac 1 R \frac {\p \Bcal_2  (X)} {\p \Gt}\right)\notag \\
& \leq  C \left| \nabla  {\Bcal_1 (X)}- b  \nabla { \Bcal_2  (X)}\right| \quad \mbox{on } l_R^+.\notag
\end{align}
In the same way, we obtain an analogous estimate on $l_R^-$. In short, we have shown that there are $R_2$ and $C>0$ independent of $R \ge R_2$ such that
\beq\label{Bcal_on_l_pm_R}
|\nabla \Bcal_1 (X)| + |\nabla \Bcal_2 (X)|  \leq  C \left| \nabla  {\Bcal_1 (X)}- b  \nabla { \Bcal_2  (X)}\right| \quad\mbox{for all } X \in l_R^+ \cup l_R^-.  
\eeq

Second, we prove \eqnref{B1B2est} on $( ( \p \GG_1 \cap B_R ) \setminus \{S_1 \})$. It follows from \eqnref{Bcal_on_l_pm_R} that
\beq\label{Bcal_on_GG_setminus_R_2}
|\nabla \Bcal_1 (X)| + |\nabla \Bcal_2 (X)|  \leq  C \left| \nabla  {\Bcal_1 (X)}- b  \nabla { \Bcal_2  (X)}\right| \quad\mbox{for all } X \in  (\p \GG_1 \cup \p \GG_2)\setminus B_{R_2}.
\eeq
By the equality \eqnref{naBcalj}, there exists a constant $C$ depending on aperture angles $\Ga_1 $ and $ \Ga_2$ such that
\beq\label{Bcal2Bcal1}
| \nabla \Bcal_2 (X) |  \le C | \nabla \Bcal_1 (X) | \quad\mbox{for all } X \in  ( \p \GG_1 \cap B_{R_2} ) \setminus \{S_1\},
\eeq
and
\beq
| \nabla \Bcal_1 (X) |  \le C | \nabla \Bcal_2 (X) | \quad\mbox{for all } X \in  ( \p \GG_2 \cap B_{R_2} ) \setminus \{S_2\}.
\eeq

Suppose $X \in ( \p \GG_1 \cap B_{R_2} ) \setminus \{S_1 \}$. Lemmas \ref{One} (i) and \eqnref{Bcal2Bcal1} show that
$$
|\nabla \Bcal_1 (X)| + | \nabla \Bcal_2 (X)| \leq  C |\nabla \Bcal_1 (X)| =  -  C \p_{\nu_1} \Bcal_1 (X).
$$
Since $b>0$ and $\p_{\nu_1} \Bcal_2\left ( X\right) \ge 0$ (the second inequality in \eqnref{Gammaone}), we have
\beq
|\nabla \Bcal_1 (X)| + | \nabla \Bcal_2 (X)| \leq - C ( \p_{\nu_1} \Bcal_1 (X) - b  \p_{\nu_1} \Bcal_2 (X)) \le C |\nabla \Bcal_1 (X) - b  \nabla \Bcal_2 (X)|. \notag
\eeq
Similarly, one can show that \eqnref{B1B2est} holds for $X \in ( \p \GG_2 \cap B_{R_2} ) \setminus \{S_2 \}$. These together with \eqref{Bcal_on_l_pm_R} and \eqref{Bcal_on_GG_setminus_R_2} complete the proof that \eqnref{B1B2est} holds on $\p \left( B_R \setminus (\GG_1 \cup \GG_2) \right) \setminus \{ S_1, S_2 \}$ for all $R \ge R_2$, as desired.

If we identify points in the plane with complex numbers, then $\nabla \Bcal_1 $, $\nabla \Bcal_2$ and $\nabla\Bcal_1 - b  \nabla\Bcal_2$ are conjugates of analytic functions. Since \eqnref{B1B2est} holds on $\p \left( B_R \setminus (\GG_1 \cup \GG_2) \right) \setminus \{ S_1, S_2 \}$, we have
$$
\left | \frac {\nabla \Bcal_1} {\nabla  {\Bcal_1 (X)}- b  \nabla { \Bcal_2  (X)}} \right| \leq C
$$
for any $X \in \p \left( B_R \setminus \overline{\GG_1 \cup \GG_2}\right) $, where $C$ is a constant independent of $R \ge R_2$. It then follows from the maximum modulus theorem that
$$
\left | \nabla \Bcal_1 (X)  \right| \leq C \left|{\nabla  {\Bcal_1 (X)}- b  \nabla { \Bcal_2  (X)}}\right| \quad\mbox{for all } X \in B_R \setminus \overline{\GG_1 \cup \GG_2}.
$$
Analogously we can derive the above inequality with $\Bcal_1$ on the left-hand side is replaced by $\Bcal_2$. Since $C$ in the estimates is independent of $R$, we arrive at \eqnref{B1B2est}, and the proof is complete. \qed

\section{Auxiliary functions on cones and their estimates}\label{sec:cone}

In this section we construct some auxiliary functions which are used in an essential way for analysing singular behavior of $\nabla u$ where $u$ is the solution to \eqnref{main}. We first construct functions in $\Rbb^2 \setminus \overline {\GG_1 \cup \GG_2}$, and then scale them by $\Ge$ in the next section.

Let, for ease of notation,
\beq
\Pi:= \Rbb^2 \setminus \overline {\GG_1 \cup \GG_2} \quad\mbox{and}\quad \Pi^+:= \{ (x,y)\in \Pi, ~y>0 \}.
\eeq
Let $P_j$ be a point on $\p\GG_j \cap \{ (x,y): y>0 \}$ for $j=1,2$ as in the previous section,
and let $Q$ be the intersection point of two straight lines containing the line segments $\overline{P_1 S_1}$ and $\overline{P_2 S_2}$. Specifically, $Q$ is given by
\beq\label{Qdef}
Q= \left(\frac {\cot \frac {\Ga_1} 2}{ \cot \frac {\Ga_1} 2 + \cot \frac {\Ga_2} 2  }- \frac 1 2, -\frac 1 { \cot   \frac {\Ga_1} 2 + \cot   \frac {\Ga_2} 2  } \right).
\eeq
We emphasize that $Q$ depends on $\Ga_1$ and $\Ga_2$. In fact, $|Q| \to \infty$ as $\Ga_1+\Ga_2 \to 2\pi$. Under the assumption \eqnref{anglecond} of this paper, we have $|Q| \approx 1$.

\begin{figure}[h!]
\begin{center}
\epsfig{figure=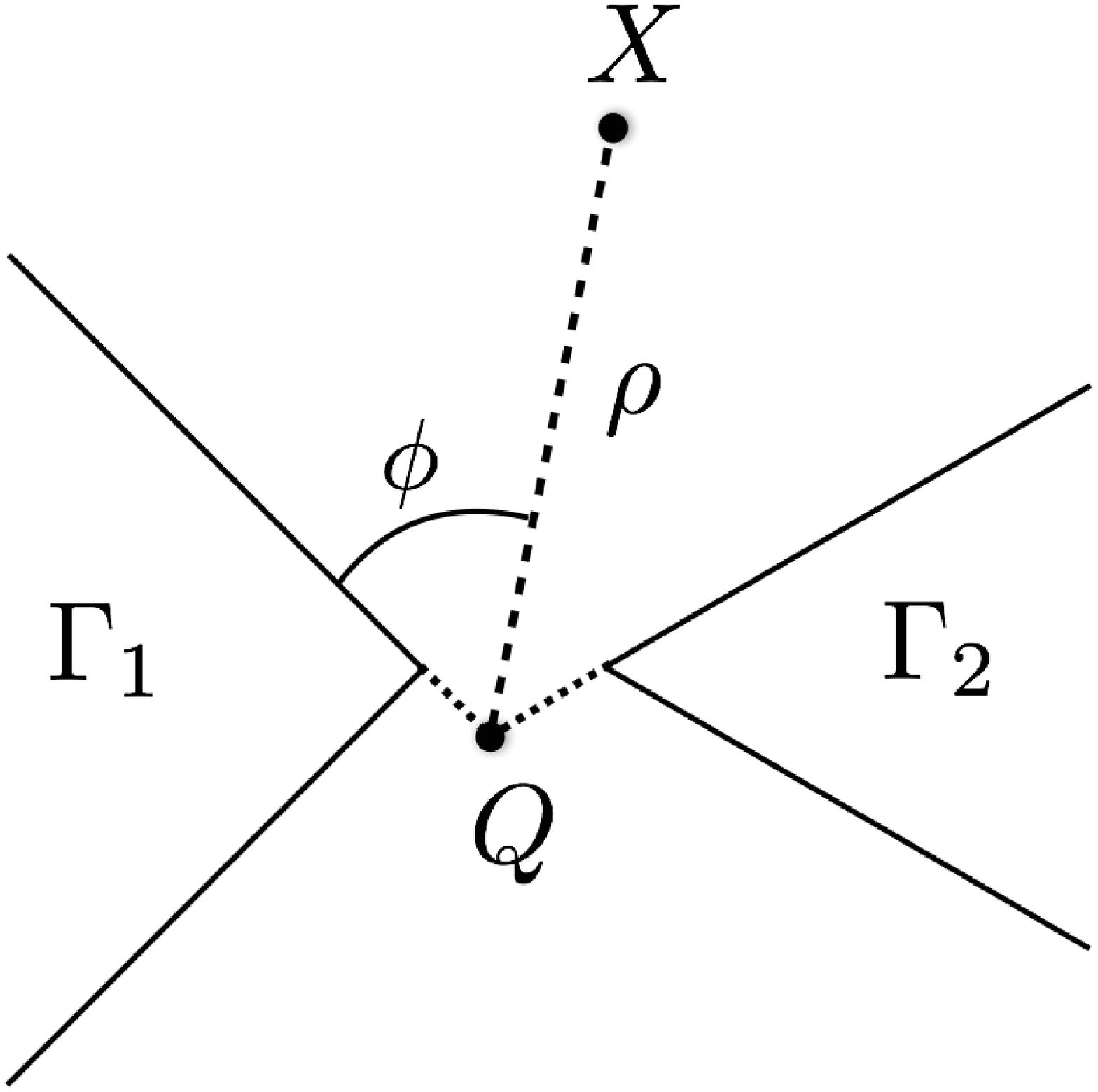,width=5cm}
\end{center}
\caption{The polar coordinates $\Gr$ and $\Gf$}\label{Fig3}
\end{figure}

Using the point $Q$, we define $\Gf(X)$ for $X \in \Pi^+$ to be the angle between $\overrightarrow {Q P_1}$ and $\overrightarrow {QX}$. Let
\beq\label{Ggdef}
\Gg := \frac {2\pi} {2\pi - (\Ga_1 +\Ga_2)}.
\eeq
Then one can easily see that
$$
0 < \Gf < \pi - \frac 1 2 (\Ga_1 + \Ga_2) = \frac{\pi}{\Gg} \quad\mbox{in } \Pi^+.
$$
We also define $\Gr=\Gr(X)$ to be the distance between $X$ and $Q$, namely,
\beq
\Gr(X):= |X-Q|.
\eeq
So, $(\Gr, \Gf)$ is the polar coordinate system with respect to $Q$ in $\Pi^+$. See Figure \ref{Fig3}.

Let $w$ be the solution to
\beq\label{wcondi}
\begin{cases}
\GD w = 0 \quad\mbox{in }  \Pi^+, \\
w= 0  \quad\mbox{on } \p\Pi^+ \setminus \left[-\frac 1 2 , \frac 1 2 \right]\times \{0\} ,  \\
\nm
\p_{y} w =  - \p_{y} \Gf \quad\mbox{on } \left[-\frac 1 2 , \frac 1 2 \right]\times \{0\} ,\\
\nm
\ds \int_{\Pi^+} |\nabla w|^2 dxdy < \infty .  \end{cases}
\eeq
We then extend $\Gf$ and $w$ to $\Pi$ as symmetric functions with respect to the $x$-axis, namely,
$$
\Gf(x,y) = \Gf(x,-y), \quad w(x,y) = w (x,-y).
$$

\begin{lem}\label{A0}
The following holds:
\begin{itemize}
\item[(i)] The function $\Gf+w$ is harmonic in $\Pi$, $\Gf+w = 0$ on $\p \GG_1$, and $\Gf+w = \pi /\Gg$ on $\p\GG_2$.
\item[(ii)] Let $\Gr_0:= \max \{ |S_1-Q|, |S_2-Q| \}$.  There is a constant $C=C(\Gr_0)$ such that
\beq\label{bound_Lemma21_2}
0<  |\nabla w (X) | \le \frac{C}{ \left| X^*- Q \right| ^{{\Gg}+1} } \quad\mbox{for all } X \in \Pi \setminus \overline{B_{2}},
\eeq
where $X^*=(x,|y|)$ for $X=(x,y)$.
\item[(iii)] For each $X \in \Pi$,
\beq\label{bound_Lemma21}
0<  (\Gf + w)(X) <  \frac{\pi}{\Gg}.
\eeq
\end{itemize}
\end{lem}

\pf
Since $\Gf+w$ is harmonic in $\Pi^+$ and $\p_{y} (\Gf+w)=0$ on $[-\frac 1 2 , \frac 1 2 ]\times \{0\}$, $\Gf+w$ is harmonic in $\Pi$. The rest of (i) follows from definitions of $\Gf$ and $w$.

To prove (ii) we first observe that the set $\{ X : \Gr>\Gr_0, \ 0< \Gf < \pi/\Gg\}$ is contained in $\Pi^+$. Because of the second condition in \eqnref{wcondi}, $w$ admits the following Fourier series expansion for $X \in \Pi^+$ satisfying $\Gr(X) > \Gr_0 $:
\beq
w (X)= \sum_{n=1}^{\infty} \frac{a_n}{\Gr^{n\Gg}} \sin \left({n\Gg \Gf}\right). \label{tilde_w_series}
\eeq
Since
\beq\label{gradw}
|\nabla w|^2 = |\p_\Gr w|^2 + \Gr^{-2} |\p_\Gf w|^2 ,
\eeq
we have
\beq\label{gradw2}
\sum_{n=1}^\infty \frac{|a_n|^2 n \Gg}{\Gr_0^{2n\Gg}} \lesssim \int_{\Gr_0}^\infty \int_{0}^{\pi/\Gg} |\nabla w|^2 \rho d\phi d\rho \le \int_{\Pi^+} |\nabla w|^2 dxdy< \infty.
\eeq
If $X \in  \Pi^+ \setminus \ol{B_2}$, then $\rho(X) \geq \Gr_0 + 1/2$, and hence
\begin{align*}
\left| \nabla w (X)\right| & \le C \sum_{n=1} ^{\infty} |a_n| \frac {n\Gg} {\Gr^{n\Gg + 1}}  \\
&\leq \frac{C}{\Gr} \left(\sum_{n=1} ^{\infty} \left|a_n\right|^2 \frac{n\Gg}{\Gr_0 ^{2 n\Gg} } \right)^{\frac 1 2 }
\left(\sum_{n=1}^{\infty}  n\Gg  \left(\frac{\Gr_0}{\Gr} \right)^{2 n\Gg}  \right)^{\frac 1 2 }
\le \frac{C}{\Gr^{\Gg+1}} ,
\end{align*}
that is,
$$
\left| \nabla w (X)\right| \le \frac{C}{ \left| X- Q \right| ^{{\Gg}+1} }.
$$
Now, \eqref{bound_Lemma21_2} follows since $w$ is symmetric with respect to the $x$-axis.

We see from \eqnref{gradw} that
$$
\p_\Gf \left(\Gf + w \right)(X) \geq 1 - |\p_\Gf w(X)| \geq 1- {\Gr(X)} | \nabla w(X)|  .
$$
Then (ii) yields
$$
\p_\Gf \left(\Gf + w \right)(X) \geq  1 - \frac{C}{\Gr^{\Gg}} > \frac 1 2,
$$
provided that $\Gr >(2C)^{1/\Gg}$. It implies that $\Gf + w$ is increasing along the arc $\{ X:\Gr=\Gr_1 \}$ where $\Gr_1$ is a constant satisfying $\Gr_1> (2C)^{1/\Gg}$. Since $\Gf + w = 0$ on $\p\GG_1$ and $\Gf + w = \pi/\Gg$ on $\p\GG_2$, we infer that \eqnref{bound_Lemma21} holds
on the arc $\{X: \Gr=\Gr_1 \}$. Since $\Gr_1$ is an arbitrary constant satisfying $\Gr_1> (2C)^{1/\Gg}$, \eqnref{bound_Lemma21} holds if $\Gr>(2C)^{1/\Gg}$.  Since $\Gf + w$ is symmetric with respect to the $x$-axis, (iii) follows from the maximal principle. This completes the proof.
\qed

\begin{lem}\label{lemma_A}
There are positive constants $a$, $b$ and $C$ such that
\beq\label{lemAest}
|\nabla \left( \Gf + w\right)(X)  - a \left(\nabla \Bcal_1 - b \nabla \Bcal_2 \right)(X) | \le C \quad \mbox{for all } X \in \Pi.
\eeq
If the aperture angles are the same, namely, $\Ga_1 = \Ga_2$, then $b=1$.
\end{lem}

\pf
We first consider the case when $X$ is near $S_1$. Since $\Gf + w=0$ on $\p \GG_1$ and the aperture angle of $\GG_1$ is $\Ga_1=2\pi-\pi/\Gb_1$, $\Gf+w$ admits the Fourier series expansion
$$
( \Gf + w ) (X)=a r_1(X)^{\Gb_1} \sin (\Gb_1 \Gt_1 (X)) + \sum_{n=2}^{\infty} a_n r_1(X)^{n\Gb_1} \sin(n\Gb_1 \Gt_1(X))
$$
in $B_1(S_1)\setminus {\overline {\GG_1}}$. Since $\sum_{n=2}^{\infty} a_n r_1(X)^{n\Gb_1} \sin(n\Gb_1 \Gt_1(X))$ converges for $r_1 <1$, its gradient converges uniformly for $r_1 \le r_0 <1$. We simply take $r_0=0.9$. Since $\Bcal_1 = r_1^{\Gb_1} \sin (\Gb_1 \Gt_1)$, we then have
\beq\label{estone}
|\nabla ( \Gf + w )(X) - a \nabla \Bcal_1 (X)| \le C
\eeq
for some constant $C$ and for all $X$ in $B_{0.9}(S_1)\setminus \overline {\GG_1}$.

Let $\GG_{1/2}$ be the arc $\{X:r_1=1/2\}$ in $B_1(S_1)\setminus {\overline {\GG_1}}$. Then, the constant $a$ is given by
$$
a = \frac{\Gb_1 2^{\Gb_1+1}}{\pi} \int_{\GG_{1/2}} \sin(\Gb_1 \Gt_1) ( \Gf + w) d\Gt_1.
$$
We then infer from \eqnref{bound_Lemma21} that
\beq
0< a < \frac{2^{\Gb_1+2}}{\Gg}.
\eeq

Applying the same argument to $\Gf + w - \pi/\Gg$ in $B_1(S_2)\setminus {\overline {\GG_2}}$, we see that there is a constant $\tilde{a}$ such that
\beq\label{esttwo}
|\nabla \left( \Gf + w \right)(X) + \tilde{a} \nabla \Bcal_2 (X) | \le C
\eeq
for all $X \in B_{0.9}(S_2) \setminus {\overline {\GG_2}}$. Since  $ 0 > \Gf + w - \pi/\Gg > -\pi/\Gg$, we obtain in the same way as before
\beq
0< \tilde{a}  < \frac{2^{\Gb_2+2}}{\Gg}.
\eeq
Let $\tilde{a}=ab$. If $\Ga_1=\Ga_2$, then $\tilde{a}=a$, and hence $b=1$.

Note that $\nabla \Bcal_2$ is bounded in $B_{0.9}(S_1)$ and $\nabla \Bcal_1$ is bounded in $B_{0.9}(S_2)$. So, it follows from \eqnref{estone} and \eqnref{esttwo} that
\beq\label{nearSl}
|\nabla \left( \Gf + w\right)(X)  - a \left(\nabla \Bcal_1 - b \nabla \Bcal_2 \right)(X) | \le C
\eeq
for all $X \in B_{0.9}(S_1) \cup B_{0.9}(S_2) \setminus \overline {\left(\GG_1 \cup \GG_2 \right)}$.

We now prove \eqnref{lemAest} away from $S_1$ and $S_2$. For doing so, set
$$
\Pi_0: = \Pi \setminus B_{0.9}(S_1) \cup B_{0.9}(S_2).
$$
According to \eqnref{bound_Lemma21}, the harmonic function $\Gf + w$ is constant on $\p \GG_1$ and $\p \GG_2$, and its range in $\Pi$ is bounded by $0$ and $\pi/\Gg$. So, for any $X_0 \in \Pi_0$, $\Gf + w$ can be extended locally as a harmonic function in $B_{1/2} (X_0)$ so that the extended function is bounded below and above by $-\pi/\Gg$ and $2\pi/\Gg$, respectively. Then, the standard gradient estimate for harmonic functions (with $R=1/2$) yields
\beq\label{standard}
\left | \nabla  \left(\Gf +  w \right) (X_0) \right| \leq \frac{2}{R} \| \Gf +  w \|_{L^\infty(B_R(X_0))} \le 2 \left( \frac 1 2 \right)^{-1}  \frac{2\pi}{\Gg}, \quad X_0 \in \Pi_0 .
\eeq
Meanwhile, since $\Pi_0$ is at distance larger than $1/2$ from $S_1$ and $S_2$, we have
$$
\| \nabla \Bcal_1 \|_{L^{\infty}(\Pi_0)} +  \| \nabla \Bcal_2 \|_{L^{\infty}(\Pi_0)} \le C.
$$
Therefore, we have for $X \in \Pi_0$
\begin{align*}
& |\nabla \left( \Gf + w\right)(X)  - a \left(\nabla \Bcal_1 - b \nabla \Bcal_2 \right)(X) | \\
& \le C(|\nabla \left( \Gf + w\right)(X)| + |\nabla \Bcal_1(X)| + |\nabla \Bcal_2 (X) |) \le C.
\end{align*}
This estimate together with \eqnref{nearSl} yield the desired estimate \eqnref{lemAest}.
\qed

\section{Singular functions and a representation of the solution}\label{sec:singular}

In this section we introduce some auxiliary functions on $\Rbb^2 \setminus \overline {(\GO_{1} \cup \GO_{2})}$, which capture the singular behavior of the gradient of the solution to \eqnref{main}.

Let $q$ be the solution to
\beq\label{qeqn}
\begin{cases}
\GD q = 0 \quad\mbox{in }  \Rbb^2 \setminus \overline {(\GO_{1} \cup \GO_{2})}, \\
\ds q = {d}_j  \quad\mbox{on } \p \GO_{j}, \ \ j=1,2, \\
\ds \int_{\p \GO_{1}} \p_{\nu} q ds = -\ds \int_{\p \GO_{2}} \p_{\nu} q ds =-1, \\
\ds q (X)  = O(|X|^{-1}) \quad\mbox{as } |X| \rightarrow \infty  .
\end{cases}
\eeq
Here, $d_1$ and $d_2$ are constants determined by the third condition in \eqnref{qeqn} and depend on $\Ge$. A proof of the existence and uniqueness of the solution to \eqnref{qeqn} can be found in \cite{ACKLY-ARMA-13}. It is worth mentioning that
\beq
q|_{\p \GO_1} = d_1 < 0 < d_2 = q |_{\p \GO_2} ,
\eeq
which can be seen using Hopf's lemma.

The singular function of a similar type was first introduced in \cite{Y} and used to estimate the gradient blow-up when $\GO_1$ and $\GO_2$ are planar domains with smooth boundaries. Recently, it has been used in an essential way to characterize the gradient blow-up when the inclusions have smooth boundaries \cite{ACKLY-ARMA-13, KLY-MA-15}. If domains have smooth boundaries, then we approximate each of them with the osculating disk. Since the singular function is explicit when $\GO_j$ are disks, we can obtain an explicit approximation of $q$.

However, since the domains under consideration in this paper have corners, such a method do not apply. In this section, we introduce some intermediate auxiliary functions which yield a good approximation of $q$.

The first intermediate function, $\Gn$, is defined to be the solution to
\beq\label{etaeqn}
\begin{cases}
\GD \Gn = 0 \quad\mbox{in }  \Rbb^2 \setminus \overline {\left(\GO_{1} \cup \GO_{2}\right)}, \\
\Gn = 0  \quad \mbox{on }\p \GO_{1},  \\
\Gn = 1    \quad \mbox{on } \p \GO_{2}  ,\\
\nm
\ds \int_{\Rbb^2 \setminus \overline {\GO_{1}\cup \GO_2}} |\nabla \Gn|^2 dxdy < \infty .  \end{cases}
\eeq
The connection between $q$ and $\Gn$ can be seen in the relations
\beq\label{qandeta}
q = (q |_{\p\GO_2} - q|_{\p\GO_1}) \Gn + q|_{\p \GO_1}
\eeq
and
\beq\label{q_and_Q}
\nabla q =  (q |_{\p\GO_2} - q|_{\p\GO_1}) \nabla \Gn.
\eeq
We emphasize that \eqnref{qandeta} shows that the solution $\eta$ to \eqnref{etaeqn} uniquely exists. Moreover, it shows that 
\beq\label{etainf}
\eta(X) \to \text{a constant}, \quad \text{as } |X| \to \infty.
\eeq

Since $\GO_j$ are symmetric with respect to the $x$-axis, so is $\Gn$, and hence $\p_\nu \Gn=0$ on the line segment between two vertices $V_1$ and $V_2$. So, in view of (i) and (iii) of Lemma \ref{A0}, it is natural to compare $\Gn$ with $\frac{\Gg}{\pi} \left( \Gf + w\right)(\Ge^{-1} X)$. For such a comparison we obtain the following lemma. We mention that the dilation $\Ge^{-1}X$ maps $X \in B_{\Gd} \setminus (\GO_1 \cup \GO_2)$ into $\Rbb^2 \setminus (\GG_1 \cup \GG_2)$. With $L_\Ge$ and $R_\Ge$ given in \eqnref{LGe}, the translation $X-L_\Ge$ maps $B_{\Gd} \cap \p\GO_1$ into $\p\GG_1$ and $X-R_\Ge$ maps $B_{\Gd} \cap \p\GO_2$ into $\p\GG_2$.

\begin{lem}\label{lemma_B_1}
Define $v(X)$ for $X \in B_{\Gd} \setminus \overline{\GO_1 \cup \GO_2}$ by
\beq\label{decomp_Q}
\Gn (X)=  \frac{\Gg}{\pi} \left( \Gf + w\right)(\Ge^{-1} X) + v (X).
\eeq
There is $\Gd_1<\Gd$ independent of $\Ge$ such that
\beq\label{nablavest}
| \nabla v (X)| \lesssim | \nabla \Bcal_1 (X-L_\Ge) | + | \nabla \Bcal_2 (X-R_\Ge) |, \quad X \in  B_{\Gd_1} \setminus \overline {\GO_1 \cup \GO_2}.
\eeq
\end{lem}

\pf
We first prove that there are constants $\Gd_1< \Gd$, $C_1$ and $C_2$ independent of $\Ge$ such that
\beq\label{bdryest}
|\nabla  v (X)| \leq C_1 \left| \nabla \Bcal_1 (X-L_\Ge) - \nabla \Bcal_2 (X - R_\Ge) + C_2 \left(x + \frac 1 2,y\right) \right |
\eeq
for all $X \in \p (B_{\Gd_1} \setminus (\GO_1 \cup \GO_2))$. Furthermore, the right-hand side of the above does not vanish.

To prove \eqnref{bdryest} on $\p\GO_1$ and $\p\GO_2$ near $O$, let $v_{1+}$ and $v_{1-}$ be the solutions to
\beq
\begin{cases}
\ds \GD v_{1+}=\GD v_{1-} =  0  \quad &\mbox{in } B_{\Gd} \setminus \overline{\GO_1}, \\
\ds v_{1+}=v_{1-} = 0   \quad &\mbox{on } \p\GO_1 \cap B_{\Gd}, \\
\ds v_{1+}=- v_{1-} = 4   \quad &\mbox{on }  \p  B_{\Gd} \setminus \overline{\GO_1},
\end{cases}
\eeq
and $v_{2+}$ and $v_{2-}$ be the solutions to
\beq
\begin{cases}
\ds \GD v_{2+}=\GD v_{2-} =  0  \quad &\mbox{in } B_{\Gd} \setminus  \overline{\GO_2},\\
\ds v_{2+}=v_{2-} = 0   \quad &\mbox{on } \p\GO_2 \cap B_{\Gd} , \\
\ds v_{2+}=- v_{2-} = 4   \quad &\mbox{on }  \p B_{\Gd} \setminus \overline{\GO_2}.
\end{cases}
\eeq
Note that $v_{1+} - v$ and $v_{1-} - v$ respectively attain the minimum and the maximum on $\p\GO_1 \cap B_{\Gd}$. So, by Hopf's lemma, we have
\beq\label{hopf1}
\p_{\nu} v_{1+}   < \p_{\nu} v < \p_{\nu} v_{1-} \quad\mbox{on } \p\GO_1 \cap B_{\Gd} .
\eeq
In the same way, we also have
\beq\label{hopf2}
\p_{\nu} v_{2+}   < \p_{\nu} v < \p_{\nu} v_{2-} \quad\mbox{on } \p\GO_2 \cap B_{\Gd} .
\eeq

Let $\Gd_2$ be such that $\ol{B_{\Gd_2}(V_1)} \setminus \ol{\GO_1} \subset B_{\Gd} \setminus \overline{\GO_1}$ and $\ol{B_{\Gd_2}(V_2)} \setminus \ol{\GO_2} \subset B_{\Gd} \setminus \overline{\GO_2}$.
Since $v_{1+}= 0$ on $\p\GO_1 \cap B_{\Gd}$, there are positive constants $a$ and $C_3$, and a function $\psi$ such that $v_{1+}$ can be expressed as
\beq\label{vone+1}
v_{1+}(X) = a \Bcal_1(X-L_\Ge) + \psi(X)
\eeq
and
\beq\label{vone+2}
|\nabla \psi(X)| \le C_3 \quad\mbox{for all } X \in \ol{B_{\Gd_2}(V_1)} \setminus \GO_1.
\eeq
In fact, $\psi$ can be expressed as
$$
\psi(X)= \sum_{n=2}^\infty a_n r_1(X-L_\Ge)^{n\Gb_1} \sin\big( n\Gb_1 \Gt_1(X-L_\Ge))
$$
for some constants $a_n$, and the Fourier series converges for $X$ satisfying $|X-L_\Ge| < \Gd_2+s$ for some $s>0$. So we achieve \eqnref{vone+2}.

Note that there is a constant $C_4>0$ such that
$$
| \nabla \Bcal_1(X-L_\Ge)| \ge C_4 \quad\mbox{for all } X \in \ol{B_{\Gd_2}(V_1)} \setminus \GO_1.
$$
So, we infer from \eqnref{vone+1} and \eqnref{vone+2} that
\begin{align*}
|\nabla v_{1+}(X)| &\le a |\nabla \Bcal_1(X-L_\Ge)| + |\nabla \psi(X)| \\
&\le (a+ C_3 C_4^{-1}) |\nabla \Bcal_1(X-L_\Ge)| \lesssim | \nabla \Bcal_1(X-L_\Ge)|
\end{align*}
for all $X \in \ol{B_{\Gd_2}(V_1)} \setminus \GO_1$.
In the same way we show that
$$
|\nabla v_{1-}(X)| \lesssim | \nabla \Bcal_1(X-L_\Ge)| \quad\mbox{for all } X \in \ol{B_{\Gd_2}(V_1)} \setminus \GO_1.
$$
It then follows from \eqnref{hopf1} that there exists a constant, say $C_5$, such that
\beq\label{1stcase_aaaa}
|\nabla  v (X)| \leq C_5 \left| \nabla \Bcal_1\left ( X-L_\Ge \right)\right |
\eeq
for all $X \in \p { \GO_1} \cap B_{\Gd_2} (V_1)$.

Since $\Bcal_1\left ( X-L_\Ge \right)=0$ for $X \in \p { \GO_1} \cap B_{\Gd_2} (V_1)$, $\nabla \Bcal_1$ in \eqnref{1stcase_aaaa} can be replaced with the normal derivative $\nu_1 \cdot \nabla \Bcal_1$. That is, we have 
\beq\label{1stcase_aaaa2}
|\nabla  v (X)| \leq C_5 \left| \p_{\nu_1} \Bcal_1\left ( X-L_\Ge \right)\right |
\eeq

Moreover, a small variation of arguments to prove \eqnref{Gammaone} shows that
\beq
\p_{\nu_1} \Bcal_1 (X-L_\Ge) <0 \quad\mbox{and}\quad -\p_{\nu_1}\Bcal_2(X-R_\Ge) \leq 0
\eeq
for all $X \in \p\GO_1 \cap B_{\Gd_2} (V_1)$. Furthermore, one can easily see that
\beq
\nu_1 \cdot \left(x + 1/2 ,y\right) <0
\eeq
for all $X \in \p\GO_1 \cap B_{\Gd_2} (V_1)$. So, we infer from \eqnref{1stcase_aaaa2} that
\begin{align}
|\nabla  v (X)| & \leq C_5 \left| \p_{\nu_1} \Bcal_1 (X-L_\Ge) -  \p_{\nu_1} \Bcal_2 (X - R_\Ge) + C \nu_1 \cdot \left(x + 1/2, y\right) \right | \nonumber \\
&\leq C_5 \left| \nu_1 \cdot \left( \nabla\Bcal_1 (X-L_\Ge) -  \nabla \Bcal_2 (X - R_\Ge) + C \left(x + 1/2, y\right) \right) \right | \nonumber \\
& \leq C_5 \left| \nabla \Bcal_1 (X-L_\Ge) - \nabla \Bcal_2 (X - R_\Ge) + C \left(x + 1/2, y\right) \right | \label{bdryest1}
\end{align}
for all $X \in \p {\GO_1} \cap B_{\Gd_2} (V_1)$.
We emphasize that \eqnref{bdryest1} holds for any positive constant $C$ and its right-hand side is non-vanishing.

Similarly to \eqnref{1stcase_aaaa}, we obtain
\beq\label{2ndcase_aaaa}
|\nabla  v (X)| \leq C_6 \left| \nabla \Bcal_2\left ( X-R_\Ge \right)\right |   \quad\mbox {for all }  X \in \p { \GO_2} \cap B_{\Gd_2} (V_2).
\eeq
It is worth mentioning that here we use \eqnref{hopf2}, not \eqnref{hopf1}. And, in a similar way, we can show that \eqnref{bdryest1} holds with $C_5$ replaced with some other constant, say $C_7$, for all $X \in \p {\GO_2} \cap B_{\Gd_2}(V_2)$ and for any positive constant $C$.

Choose $\Gd_3$ so that $B_{\Gd_3} \subset B_{\Gd_1} (V_1) \cap B_{\Gd_2} (V_2)$. Then, we arrive at
\beq\label{bdryest4}
|\nabla  v (X)| \leq C_8 \left| \nabla \Bcal_1 (X-L_\Ge) - \nabla \Bcal_2 (X - R_\Ge) + C \left(x + 1/2, y\right) \right |
\eeq
for all $X \in (\p\GO_1 \cup \p\GO_2) \cap B_{\Gd_3}$ and for any positive constant $C$.

We now prove \eqnref{bdryest} on $\p B_{\Gd_1} \setminus (\GO_1 \cup \GO_2)$ if $\Gd_1<\Gd_3$.
By the maximum principle, we have $|\Gn| \le 1$. We also have $|\frac{\Gg}{\pi}(\Gf + w)| \le 1$ from \eqnref{bound_Lemma21}. So it follows from the definition \eqnref{decomp_Q} of $v$ that
\beq
|v(X)| \le 2, \quad X \in B_{\Gd_3} \setminus \overline{\GO_1 \cup \GO_2}.
\eeq
The function $v$ takes constant values on $B_{\Gd_3} \cap \p\GO_1$ and $B_{\Gd_3} \cap \p\GO_2$, and
\beq\label{vvalue}
v|_{B_{\Gd_3} \cap \p\GO_1}= v|_{B_{\Gd_3} \cap \p\GO_2}.
\eeq
It is helpful to mention here that $v$ actually takes value $0$ on $B_{\Gd_3} \cap \p\GO_1$ and $B_{\Gd_3} \cap \p\GO_2$. But we only use the weaker property \eqnref{vvalue}.
We then infer from the standard gradient estimate as used in \eqnref{standard} that if $\Gd_1<\Gd_3$, then there is a constant $C_9$ (independent of $\Ge$) such that
\beq\label{3rdcase_aaaa}
|\nabla  v (X)| \leq C_9 \quad\mbox{for all } X \in \p B_{\Gd_1} \setminus (\GO_1 \cup \GO_2).
\eeq

Since $|\nabla \Bcal_1 (X-L_\Ge) - \nabla \Bcal_2 (X - R_\Ge)|$ is bounded for all $X \in \p B_{\Gd_1} \setminus (\GO_1 \cup \GO_2)$, we may choose the constant $C$ appearing in \eqnref{bdryest4} large enough (and denote it by $C_2$) so that
$$
C_8 \left| \nabla \Bcal_1 (X-L_\Ge) - \nabla \Bcal_2 (X - R_\Ge) + C_2 \left(x + 1/2, y\right) \right | \ge C_9
$$
for all $X \in \p B_{\Gd_1} \setminus (\GO_1 \cup \GO_2)$, where $C_8$ and $C_9$ are constants appearing in \eqnref{bdryest4} and \eqnref{3rdcase_aaaa}, respectively. Then, we have
\beq\label{bdryest3}
|\nabla  v (X)| \leq C_8 \left| \nabla \Bcal_1 (X-L_\Ge) - \nabla \Bcal_2 (X - R_\Ge) + C_2 \left(x + 1/2, y\right) \right |
\eeq
for all $X \in \p B_{\Gd_1} \setminus (\GO_1 \cup \GO_2)$. Since \eqnref{bdryest4} holds for any $C>0$ on $(\p\GO_1 \cup \p\GO_2) \setminus B_{\Gd_1}$, \eqnref{bdryest} is proved with $C_1=C_8$.

If we identify points in the plane with complex numbers, then $\nabla v$ and $\nabla \Bcal_1 (X-L_\Ge) - \nabla \Bcal_2 (X - R_\Ge) + C_2 \left(x + 1/2,y\right)$ are conjugates of analytic functions. Furthermore, the latter function does not vanish in $B_{\Gd_1} \setminus (\GO_1 \cup \GO_2)$ due to the maximum modulus theorem. So, we see that
$$
\left|\frac {\nabla  v (X) } {\nabla \Bcal_1 ( X - L_\Ge) - \nabla \Bcal_2 ( X - R_\Ge ) + C_2 \left(x+1/2,y\right)}\right| \leq C, \quad X \in B_{\Gd_1} \setminus (\GO_1 \cup \GO_2)
$$
for some constant $C$. It then follows from the triangular inequality that
$$
|\nabla  v (X)| \lesssim |\nabla \Bcal_1 ( X - L_\Ge)| + |\nabla \Bcal_2 ( X - R_\Ge )| + \left| \left( x+1/2,y \right) \right|.
$$
Since $|\nabla \Bcal_1 ( X - L_\Ge)| + |\nabla \Bcal_2 ( X - R_\Ge )| \gtrsim 1$ while $|(x+1/2,y)| \lesssim 1$ for all $X \in B_{\Gd_1} \setminus (\GO_1 \cup \GO_2)$, \eqnref{nablavest} follows and the proof is complete. \qed

\medskip

Define the constant $c(u)$ by
\beq\label{bdef}
c(u):= \frac {u |_{\p\GO_2} - u |_{\p\GO_1}}{q |_{\p\GO_2} - q|_{\p\GO_1}},
\eeq
and define the function $\Gs$ by
\beq\label{decomp_u=bq+r}
u = c(u) q + \Gs.
\eeq
Then $\Gs$ is harmonic in $\Rbb^2 \setminus \overline {(\GO_{1} \cup \GO_{2})}$ having a constant value on $\p\GO_1 \cup \p \GO_2$ such that
\beq\label{Gsbdry}
\Gs |_{\p \GO_1} = \Gs|_{\p \GO_2} .
\eeq

We prove the following lemma.

\begin{lem}\label{lemma_B}
Let $\GO$ be a bounded set containing $\ol{B_{\Gd}}$ and $\overline{\GO_1 \cup \GO_2}$. There is $\Gd_1$ independent of $\Ge$ such that
\beq\label{naGs}
| \nabla \Gs (X) | \lesssim \| h \|_{L^{\infty}(\GO)} \left(\left| \nabla \Bcal_1\left ( X-L_\Ge \right)\right | + \left| \nabla \Bcal_2\left ( X-R_\Ge \right) \right |\right)
\eeq
for all $X \in  B_{\Gd_1} \setminus \overline {\GO_1 \cup \GO_2}$.
\end{lem}
\pf
Thanks to \eqnref{Gsbdry}, \eqnref{naGs} can be proved in the same way as the proof of Lemma \ref{lemma_B_1}, once we prove that
\beq\label{normGs}
\norm{\Gs}_{L^{\infty} (B_{\Gd} \setminus \overline{\GO_1 \cup \GO_2})} \lesssim \| h \|_{L^{\infty}(\GO)} .
\eeq

To prove \eqnref{normGs}, let
$$
\Gn_0 = \lim_{|X| \to \infty} \Gn(X),
$$
where the limit exists as seen at \eqnref{etainf}. Then we have $0< \Gn_0 < 1$, and
\beq\label{qeta}
q(X) = \left(q|_{\p \GO_2} - q|_{\p \GO_1}\right) \left(\Gn(X) - \Gn_0\right), \quad X \in \Rbb^2 \setminus \overline{\GO_1 \cup \GO_2}.
\eeq
In fact, we denote the right-hand side of the above by $f$, then $f$ is harmonic in $\Rbb^2 \setminus \overline{\GO_1 \cup \GO_2}$, $f(X) \to 0$ as $|X| \to \infty$, $f$ is constant on $\p\GO_j$ for $j=1,2$, and $f|_{\p\GO_2} - f|_{\p\GO_1}=q|_{\p\GO_2} - q|_{\p\GO_1}$. So, \eqnref{qeta} holds.
We then have
\beq\label{eq_1emma_B}
\Gs(X) = \left ( u |_{\p \GO_1} -  u |_{\p \GO_2}\right) \left(\Gn(X) - \Gn_0\right) + u(X), \quad X \in \Rbb^2 \setminus \overline{\GO_1 \cup \GO_2}.
\eeq

It was proven in \cite{Y} that
\beq\label{the_presentation_potential_diff_u}
u |_{\p \GO_2} -  u |_{\p \GO_1} = \int_{\p \GO_1 \cup \p \GO_2} h \p_{\nu} q ds.
\eeq
Thus, we have
$$
\left|  u |_{\p \GO_2} -  u |_{\p \GO_1} \right|\leq  \norm{h}_{L^{\infty}(\GO)} \int_{\p\GO_1 \cup \p\GO_2} |\p_{\nu} q| ds .
$$
Since $\p_{\nu} q$ is either positive or negative on $\p\GO_j$ for $j=1,2$, we have
\beq\label{900}
\left|  u |_{\p \GO_2} -  u |_{\p \GO_1} \right| \leq \norm{h}_{L^{\infty}(\GO)} \left(\left | \int_{\p \GO_1} \p_{\nu} q ds \right | + \left| \int_{ \p \GO_2} \p_{\nu} q ds \right| \right) \leq 2  \norm{ h}_{L^{\infty}(\GO)} .
\eeq
Since $0< \Gn \leq 1$ and $0< \Gn_0 < 1$, it follows that
\beq\label{1000}
\norm {\left ( u |_{\p \GO_2} -  u |_{\p \GO_1}\right) \left(\Gn - \Gn_0\right) }_{L^{\infty} (B_{\Gd} \setminus \overline{\GO_1 \cup \GO_2})} \leq 4  \norm{h} _{L^{\infty} (\GO)}.
\eeq

Since $u(X)-h(X) \to 0$ as $|X| \to \infty$, we have
\begin{align*}
\norm{u-h}_{L^{\infty} (B_{\Gd} \setminus \overline{\GO_1 \cup \GO_2})} & \le \norm{u-h}_{L^{\infty} ({\p \GO_1 \cup \p \GO_2})} \\
& \le \max_{\p \GO_1 \cup \p \GO_2} (u-h) - \min_{\p \GO_1 \cup \p \GO_2} (u-h) \\
& \le \left|  u |_{\p \GO_2} - u |_{\p \GO_1} \right|  + 2\norm{h}_{L^{\infty}(\GO)}.
\end{align*}
We then obtain from \eqnref{900}
\beq\label{1100}
\norm{u}_{L^{\infty}(B_{\Gd} \setminus \overline{\GO_1 \cup \GO_2})} \leq \norm {u-h}_{L^{\infty} (B_{\Gd} \setminus \overline{\GO_1 \cup \GO_2})} + \norm{h}_{L^{\infty}(B_{\Gd} \setminus \overline{\GO_1 \cup \GO_2})} \lesssim \norm{h}_{L^{\infty}(\GO)}.
\eeq
Then \eqnref{normGs} follows from \eqnref{eq_1emma_B}, \eqnref{1000} and \eqnref{1100}.
\qed

\begin{lem}\label{lemma_C}
There is $\Ge_0 >0$ such that
\beq
q|_{\p \GO_2} - q|_{\p \GO_1} \simeq \frac 1 {|\log \Ge|}
\eeq
for all $\Ge \le \Ge_0$.
\end{lem}

\pf
We first observe from \eqnref{q_and_Q} that
$$
(q |_{\p \GO_2} - q|_{\p \GO_1}) \int_{\p \GO_2}  \p_{\nu} \Gn ds = \int_{\p \GO_2} \p_{\nu} q ds = 1.
$$
So, it suffices to prove
\beq\label{main_lem_2-5}
\int_{\p \GO_2}  \p_{\nu} \Gn ds \simeq  {|\log \Ge|}.
\eeq

To prove \eqnref{main_lem_2-5}, we write
\beq\label{1st+2nd_lem_2-5}
\int_{\p \GO_2} \p_{\nu} \Gn ds = \int_{\p \GO_2 \setminus B_{\Gd_1}} \p_{\nu} \Gn  ds  + \int_{\p \GO_2 \cap B_{\Gd_1}} \p_{\nu} \Gn ds =: I+II,
\eeq
where $\Gd_1$ is the number appearing in Lemmas \ref{lemma_B_1} and \ref{lemma_B}. Since $\Gn$ is constant on $\p \GO_2$, $0\le \Gn \le 1$, and the set $\p \GO_2 \setminus B_{\Gd_1}$ is at some distance from the vertex $V_2$, there is $r_0$ independent of $\Ge$ such that $\eta$ is extended by reflection as a harmonic function in $B_{r_0}(X)$ for each $X \in \p \GO_2 \setminus B_{\Gd_1}$ by a conformal mapping and the extended function, which is still denoted by $\Gn$, satisfies $-1 \le \Gn \le 2$. So, we have
$$
\left| \nabla \Gn(X) \right| \lesssim \frac{2}{r_0} \norm{\Gn}_{L^{\infty}(B_{r_0}(X))} \lesssim \frac{4}{r_0}
$$
for any $X \in \p\GO_2 \setminus B_{\Gd_1}$. Thus, there exists a constant $M_1$ independent of $\Ge$ and $h$ such that
\beq\label{I}
I \leq M_1.
\eeq

In view of \eqnref{decomp_Q}, we may write $II$ as
\begin{align*}
II & = \int_{\p \GO_2 \cap B_{\Gd_1}} \p_{\nu} v ds + \frac{\Gg}{\pi} \int_{\p \GO_2 \cap B_{3\Ge}} \p_{\nu} \left(  (\Gf + w) (\Ge^{-1} X) \right) ds\\
& \qquad + \frac{\Gg}{\pi} \int_{\p \GO_2 \cap (B_{\Gd_1} \setminus B_{3\Ge} )} \p_{\nu} \left( w (\Ge^{-1} X)\right)  ds
+ \frac{\Gg}{\pi} \int_{\p \GO_2 \cap (B_{\Gd_1} \setminus B_{3\Ge} )} \p_{\nu} \left( \Gf (\Ge^{-1} X)\right)  ds \\
&=: II_1 + II_2 + II_3 + II_4.
\end{align*}
By Lemma \ref{lemma_B_1}, we have
\beq\label{II1}
|II_1| \lesssim \int_{\p \GO_2 \cap B_{\Gd_1}} | \nabla \Bcal_1 (X-L_\Ge) | + | \nabla \Bcal_2 (X-R_\Ge) | \, ds \lesssim 1,
\eeq
where the last inequality holds since $\Gb_j >0$. If we make a change of variables $\Ge^{-1} X \to X$, it follows from Lemma \ref{lemma_A} that
\begin{align}
|II_2| &\le \frac{\Gg}{\pi} \int_{(\p \GG_2) \cap B_{3}} |\p_{\nu} (\Gf + w ) (X)| ds \nonumber \\
&\lesssim \int_{(\p \GG_2) \cap B_{3}} (|(\nabla \Bcal_1 - b \nabla \Bcal_2)(X)| + 1) ds \lesssim 1. \label{II2}
\end{align}

The same change of variables yields that
$$
|II_3| \lesssim \int_{\p \GG_2 \setminus B_{3}} \left | \p_{\nu} w (X) \right|  ds \le 2 \int_{\p \GG_2^+ \setminus B_{5/2}(S_2)} \left | \p_{\nu} w (X) \right|  ds,
$$
where $\p \GG_2^+ = \{ X \in \p \GG_2: y>0 \}$. Let $\Gr_0:= \max \{ |S_1-Q|, |S_2-Q| \}$ as in Lemma \ref{A0} (ii). Since $|S_1-S_2|=1$, the triangular inequality shows that
$$
\min \{ |S_1-Q|, |S_2-Q| \} \ge \Gr_0-1.
$$
So, we have from \eqnref{gradw} that
\begin{align*}
|II_3| &\lesssim \int_{\p \GG_2^+ \setminus B_{3/2+\Gr_0}(Q)} \left | \p_{\nu} w (X) \right|  ds \\
& \lesssim \sum_{n=1}^{\infty} |a_n| n\Gg \int_{3/2+\Gr_0}^\infty \frac{1}{\Gr^{(n\Gg + 1)}} d\Gr \lesssim \sum_{n=1}^{\infty} \frac{|a_n|}{(3/2+\Gr_0)^{n\Gg}} .
\end{align*}
It then follows from the Cauchy-Schwartz inequality and \eqnref{gradw2} that
\beq\label{II3}
|II_3| \lesssim \left( \sum_{n=1}^\infty \frac{|a_n|^2 n \Gg}{\Gr_0^{2n\Gg}} \right)^{1/2} \left( \sum_{n=1}^\infty \frac{\Gr_0^{2n\Gg}}{(3/2+\Gr_0)^{2n\Gg}} \right)^{1/2} < \infty.
\eeq

We apply the same change of variables $\Ge^{-1} X \to X$ to see that
$$
II_4 \simeq \Gg \int_{\p \GG_2 \cap (B_{1/(2\Ge)} \setminus B_{3})} |\p_{\nu} \Gf (X)| ds.
$$
So, we have
\beq\label{II4}
II_4 \simeq \Gg \int_{3}^{\frac{1}{2\Ge}} |\p_{\nu} \Gf (X)|d\Gr =  \Gg \int_{3}^{\frac{1}{2\Ge}} \Gr^{-1} |\p_\Gf \Gf(X)| d\Gr \simeq |\log \Ge| .
\eeq
Then \eqref{main_lem_2-5} follows from \eqnref{I}-\eqnref{II4}. This completes the proof.
\qed

\begin{lem}\label{lemma_for_b}
The constant $c(u)$ defined by \eqref{bdef} satisfies
\beq\label{cuest}
|c(u)| \lesssim \norm{h}_{L^{\infty}(\GO)} .
\eeq
Here, $\GO$ is the domain appearing in Lemma \ref{lemma_B}. 
\end{lem}
\pf
From \eqnref{q_and_Q}, \eqref{bdef}, and \eqref{the_presentation_potential_diff_u}, we have
\beq
c(u) =  \frac {u |_{\p \GO_2} - u|_{\p \GO_1}}{q |_{\p \GO_2} - q|_{\p \GO_1} } = \int_{\p \GO_1 \cup \p \GO_2} h \p_{\nu} \Gn ds = \int_{\p \GO_1 \cup \p \GO_2} (h(X) - h(O)) \p_{\nu} \Gn ds .
\eeq
We only deal with the integral over $\p\GO_2$ to prove that it is bounded independently of $\Ge$. The integral over $\p\GO_1$ can be dealt with in the same way.

Using \eqnref{decomp_Q}, we write
\begin{align*}
& \int_{\p \GO_2} (h(X) - h(O)) \p_{\nu} \Gn ds \\
&= \int_{\p\GO_2 \setminus B_{\Gd_1}} (h(X) - h(O)) \p_{\nu} \Gn ds +  \int_{\p \GO_2 \cap B_{\Gd_1}} (h(X) - h(O)) \p_{\nu} v ds \\
&\quad + \frac{\Gg}{\pi} \int_{\p \GO_2 \cap B_{3\Ge}} (h(X) - h(O)) \p_{\nu} \left( \left( \Gf+w \right) \left(\Ge^{-1} X \right)  \right) ds \\
&\quad + \frac{\Gg}{\pi} \int_{\p \GO_2 \cap (B_{\Gd_1} \setminus B_{3\Ge} )} (h(X) - h(O)) \p_{\nu} \left( w (\Ge^{-1} X)\right)  ds \\
&\quad + \frac{\Gg}{\pi} \int_{(\p \GO_2) \cap (B_{\Gd_1} \setminus B_{3\Ge} )}  (h(X) - h(O)) \p_{\nu} ( \Gf (\Ge^{-1} X) ) ds  =: \sum_{j=1}^5 I_j.
\end{align*}

Using the same arguments for the proof of Lemma \ref{lemma_C}, one can show that
$$
|I_j| \lesssim \norm{h}_{L^{\infty} (\GO)}, \quad j=1, \cdots, 4.
$$

To estimate $I_5$, let $\p \GO_2^{\pm} = \{ X \in \p\GO_2 ~:~\pm y \geq 0 \}$.
We then have
\begin{align*}
&\left| \int_{\p \GO_2^+ \cap (B_{\Gd_1} \setminus B_{3\Ge})} (h(X) - h(O)) \p_{\nu} (\Gf (\Ge^{-1} X) ) ds \right| \\
&= \left|\int_{\p \GG_2^+ \cap (B_{\Gd_1/\Ge} \setminus B_{3} )} ( h(\Ge X) - h(O)) \p_{\nu} \Gf (X) ds \right| \\
&\leq \norm{\nabla h} _{L^{\infty} (B_{\delta_1})} \left| \int_{\p \GG_2^+ \cap (B_{\Gd_1/\Ge} \setminus B_{3/2(V_2)} )} \Ge |X|  \Gr^{-1} \p_{\Gf} \Gf(X) d\Gr \right| \\
&\lesssim \norm{\nabla h}_{L^{\infty} (B_{\delta_1})} \left| \int_{\frac 3 2}^{\frac{\Gd_1}{\Ge}}  \Ge d\Gr \right| \lesssim \norm{\nabla h}_{L^{\infty} (B_{\delta_1})} .
\end{align*}
The integral over $\p \GO_2^- \cap (B_{\Gd_1} \setminus B_{3\Ge})$ can be dealt with in the exactly same manner. So, we have
$$
|I_5| \lesssim \norm{\nabla h}_{L^{\infty} (B_{\delta_1})} \lesssim  \norm{ h}_{L^{\infty} (\GO)}.
$$
This completes the proof.
\qed

\section{Estimates for the field enhancement}\label{sec:main}

In this section we present main estimates of the gradient $\nabla u$ of the solution $u$ to \eqnref{main}. We derive two different estimates, one for points quite close to the vertices and the other for those relatively away from the vertices.

Let us first consider the case when $X$ is relatively away from the vertices.
By \eqref{q_and_Q} and \eqref{decomp_u=bq+r}, we have
$$
\nabla u = c(u) \nabla q+ \nabla \Gs = c(u) (q |_{\p\GO_2} - q|_{\p\GO_1}) \nabla \Gn + \nabla \Gs.
$$
It then follows from \eqref{decomp_Q} that
\begin{align}
\nabla u (X) &= c(u) (q |_{\p \GO_2} - q|_{\p \GO_1}) \nabla \left( \frac{\Gg}{\pi} ( \Gf + w ) (\Ge^{-1} X) + v (X) \right) + \nabla \Gs (X) \nonumber \\
&= c(u) (q |_{\p \GO_2} - q|_{\p \GO_1}) \left[ \frac{\Gg}{\Ge\pi} ( \nabla \Gf + \nabla w ) (\Ge^{-1} X) + \nabla v (X) \right] + \nabla \Gs (X). \label{3000}
\end{align}

We now show that $(\nabla \Gf) (\Ge^{-1} X)$ is the major term in describing singular behavior of $\nabla u$. To do so, let
\beq\label{aGedef}
a_\Ge := \frac{c(u) (q |_{\p \GO_2} - q|_{\p \GO_1})\Gg |\log \Ge|}{\pi} = \frac{(u |_{\p \GO_2} - u|_{\p \GO_1})\Gg |\log \Ge|}{\pi}.
\eeq
Lemmas \ref{lemma_C} and \ref{lemma_for_b} show that
\beq\label{aGe}
|a_{\Ge}| \lesssim \norm{h}_{L^{\infty}(\GO)}.
\eeq
This is a crucial fact in determining a upper bound on $|\nabla u|$ in the sequel. A lower bound on $|a_{\Ge}|$ also determines a lower bound on $|\nabla u|$. However, $a_\Ge$ may or may not be $0$ depending on the configuration of inclusions and the background potential $h$. We show in section \ref{sec:aGe} that $a_\Ge$ is bounded below if $h(x,y)=x$, given the configuration of the bow-tie structure in the paper.

Let
\beq\label{Rcaldef}
\Rcal(X) := c(u) (q |_{\p \GO_2} - q|_{\p \GO_1})  \nabla v(X)+ \nabla \Gs(X).
\eeq
Then we obtain the following theorem.

\begin{thm}\label{1st_thm_two}
Let $u$ be the solution to \eqnref{main}. Then $\nabla u$ admits the following decomposition:
\beq\label{firstrep}
\nabla u (X)= \frac{a_\Ge}{\Ge |\log \Ge|} \left[ (\nabla\Gf)(\Ge^{-1} X) + (\nabla w ) (\Ge^{-1} X) \right] + \Rcal(X), \quad  X \in B_{\Gd_1} \setminus \overline{ \GO_1 \cup \GO_2}.
\eeq
Moreover, the following estimates hold for all $X \in B_{\Gd_1} \setminus \overline{ \GO_1 \cup \GO_2}$ with $|X| \ge 2\Ge$, where $\Gd_1$ is the constant appearing in Lemma \ref{lemma_B}:
\begin{align}
|(\nabla\Gf)(\Ge^{-1} X)| &\simeq \frac{\Ge}{|X|}, \label{2002} \\
|(\nabla w)(\Ge^{-1}X) | &\lesssim \frac{\Ge^{\Gg+1}}{|X|^{{\Gg}+1}}, \label{2001} \\
| \Rcal(X) | &\lesssim \norm{h}_{L^{\infty} (\GO)} (|X|^{\Gb_1-1} + |X|^{\Gb_2-1}). \label{2000}
\end{align}
\end{thm}
\pf
The decomposition \eqnref{firstrep} is an immediate consequence of \eqnref{3000}-\eqnref{Rcaldef}.

If $X=(x,y)\in B_{\Gd_1} \setminus \overline{ \GO_1 \cup \GO_2}$ is such that $|X| \ge 2\Ge$ and $y>0$, then 
$$
\Gr(\Ge^{-1}X)= |\Ge^{-1} X-Q| \ge \Gr_0
$$ 
where $\Gr_0$ is the number appearing in Lemma \ref{A0} (ii), and 
$$
|\Ge^{-1} X-Q| \simeq \Ge^{-1} |X|.
$$ 
So, it follows from \eqnref{bound_Lemma21_2} that
$$
| \left(\nabla w\right)(\Ge^{-1}X) | \lesssim \frac{1}{\left| \Ge^{-1} X- Q \right|^{{\Gg}+1}} \lesssim \frac{\Ge^{\Gg+1}}{|X|^{{\Gg}+1}}.
$$
Since $|\nabla \Gf| = \Gr^{-1}$, we have
$$
|(\nabla\Gf)(\Ge^{-1} X)| = \Gr(\Ge^{-1} X)^{-1} \simeq \Ge |X|^{-1}.
$$
So, we have \eqnref{2002} and \eqnref{2001} when $y>0$. Since $\Gf$ and $w$ are symmetric with respect to the $x$-axis, we obtain \eqnref{2002} and \eqnref{2001} for $y<0$ as well.

According to Lemmas \ref{lemma_B_1}, \ref{lemma_B}, \ref{lemma_C}, and \ref{lemma_for_b}, we have
\beq\label{2000-1}
|\Rcal(X) | \lesssim \norm{h}_{L^{\infty} (\GO)} \left(\left| \nabla \Bcal_1 ( X-L_\Ge ) \right | + \left| \nabla \Bcal_2 ( X-R_\Ge ) \right |\right)
\eeq
for any $X \in B_{\Gd_1} \setminus \overline{\GO_1 \cup \GO_2}$. Note that $L_\Ge + S_1=V_1$. So, it follows from \eqnref{naBcalj} that
$$
|\nabla \Bcal_1 (X-L_\Ge)| = \Gb_1 |X-V_1|^{\Gb_1-1}.
$$
Likewise, we have
$$
|\nabla \Bcal_2 (X-R_\Ge)| = \Gb_2 |X-V_2|^{\Gb_2-1}.
$$
So, we have
\beq\label{2000-2}
|\Rcal(X) | \lesssim \norm{h}_{L^{\infty} (\GO)} \left( |X-V_1|^{\Gb_1-1} + |X-V_2|^{\Gb_2-1} \right)
\eeq
for any $X \in B_{\Gd_1} \setminus \overline{\GO_1 \cup \GO_2}$. If $|X| \ge 2\Ge$, then $|X-V_j| \simeq |X|$. So, \eqnref{2000} follows from \eqnref{2000-2}. This completes the proof. \qed

\medskip

Let us now derive estimates of $|\nabla u|$. We first observe from Theorem \ref{1st_thm_two} that the following holds:
$$
|\nabla u(X)| \lesssim \frac{\norm{h}_{L^{\infty} (\GO)}}{|\log \Ge|} |X|^{-1}
$$
as long as $X \in B_{\Gd_1} \setminus \overline{ \GO_1 \cup \GO_2}$ and $2 \Ge \le |X| \le |\log\Ge|^{-1/\Gb}$, where $\Gb:= \min \{\Gb_1, \Gb_2 \}$. 

To derive the opposite inequality, we see from \eqnref{2002}-\eqnref{2000} that there are positive constants, say $C_\Gf$, $C_w$ and $C_\Rcal$, independent of $\Ge$ and $X$ such that
$$
|(\nabla\Gf)(\Ge^{-1} X)| \ge \frac{C_\Gf \Ge}{|X|}, \quad
|(\nabla w)(\Ge^{-1}X) | \le \frac{C_w \Ge^{\Gg+1}}{|X|^{{\Gg}+1}},
$$
and
$$
| \Rcal(X) | \le C_\Rcal (|X|^{\Gb_1-1} + |X|^{\Gb_2-1})
$$
as long as $X \in B_{\Gd_1} \setminus \overline{ \GO_1 \cup \GO_2}$ and $|X| \ge 2\Ge$. We then have
\begin{align*}
\big| (\nabla\Gf)(\Ge^{-1} X) - (\nabla w)(\Ge^{-1}X) \big| &\ge |(\nabla\Gf)(\Ge^{-1} X)| - |(\nabla w)(\Ge^{-1}X)| \\
&\ge \frac{C_\Gf \Ge}{|X|} - \frac{C_w \Ge^{\Gg+1}}{|X|^{{\Gg}+1}}.
\end{align*}
So, if $|X|^\Gg \ge 2 C_\Gf^{-1} C_w \Ge^\Gg$, we have
$$
\big| (\nabla\Gf)(\Ge^{-1} X) - (\nabla w)(\Ge^{-1}X) \big| \ge \frac{C_\Gf \Ge}{2|X|}.
$$
It then follows from \eqnref{firstrep} that
\begin{align*}
|\nabla u(X)| &\ge \frac{|a_\Ge|}{\Ge |\log \Ge|} \big| (\nabla\Gf)(\Ge^{-1} X) - (\nabla w)(\Ge^{-1}X) \big| - |\Rcal(X)| \\
&\ge \frac{|a_\Ge| C_\Gf}{2|\log \Ge||X|} - 2 C_\Rcal |X|^{\Gb-1}.
\end{align*}
Here we assumed $|X| \le 1$. Suppose that $|a_\Ge| \ge a_0$ for some positive constant $a_0$. In such a case, if $|X|^\Gb \le \frac{1}{8} a_0 C_\Rcal^{-1} C_\Gf |\log\Ge|^{-1}$, then we have
$$
|\nabla u(X)| \ge \frac{|a_\Ge| C_\Gf}{4|\log \Ge||X|}.
$$
In short, we obtain the following theorem.

\begin{thm}\label{cortwo}
Let $\Gb:= \min \{\Gb_1, \Gb_2 \}$. 
\begin{itemize}
\item[(i)] We have
\beq\label{|X|^{-1}2}
|\nabla u(X)| \lesssim \frac{\norm{h}_{L^{\infty} (\GO)}}{|\log \Ge|} |X|^{-1}
\eeq
for all $X \in B_{\Gd_1} \setminus \overline{ \GO_1 \cup \GO_2}$ with $2 \Ge \le |X| \le |\log\Ge|^{-1/\Gb}$.

\item[(ii)] Suppose that $|a_\Ge| \ge a_0$ for some positive constant $a_0$. There are positive constants $c_1 \ge 2$ and $c_2$ independent of $\Ge$ and $X$ such that
\beq\label{|X|^{-1}22}
\frac{1}{|\log \Ge|} |X|^{-1} \lesssim |\nabla u(X)|
\eeq
for all $X$ satisfying $c_1 \Ge \le |X| \le c_2 |\log\Ge|^{-1/\Gb}$.
\end{itemize}
\end{thm}

Theorem \ref{cortwo} reveals that the singularity of $\nabla u(X)$ is stronger than the corner singularity which is of the order $|X|^{-1+\Gb_1}$ or $|X|^{-1+\Gb_2}$.

We now look into the case when $X$ is close to vertices and obtain the following theorem for asymptotic behavior of $\nabla u(X)$.

\begin{thm}\label{1st_thm}
Let $u$ be the solution to \eqnref{main}.  Then $\nabla u$ admits the decomposition
\beq\label{naurep}
\nabla u(X) = \frac{a_{\Ge}a}{\Ge |\log \Ge|} (\nabla \Bcal_1 - b \nabla \Bcal_2) ( \Ge^{-1} X ) + \Ecal(X), \quad X \in B_{\Gd_1} \setminus \overline{\GO_1 \cup \GO_2},
\eeq
where $a$ and $b$ are positive constants independent of $\Ge >0$ and $h$. In particular, if two aperture angles are the same, namely, $\Ga_1 = \Ga_2$, then $b=1$. The error term $\Ecal$ satisfies
\beq\label{naurep2}
|\Ecal(X)| \lesssim \norm{h}_{L^{\infty} (\GO)} \left(\frac{1}{\Ge |\log \Ge|} + \sum_{j=1}^2 |X-V_j|^{\Gb_j-1}\right).
\eeq
\end{thm}

\pf
Let
\beq
\Ncal(X) := \frac{1}{a} \nabla (\Gf + w)(X) - \left(\nabla \Bcal_1 - b \nabla \Bcal_2 \right)(X).
\eeq
Then, we see from \eqnref{firstrep} that \eqnref{naurep} holds with
\beq\label{error}
\Ecal(X) = \frac{a_{\Ge}a}{\Ge |\log \Ge|} \Ncal(\Ge^{-1}X) + \Rcal(X).
\eeq

We infer from \eqnref{lemAest} that
$$
|\Ncal(\Ge^{-1}X)| \lesssim 1,
$$
which together with \eqnref{2000-2} yields \eqnref{naurep2}. This completes the proof. \qed

Since
\beq\label{4000}
|(\nabla \Bcal_j) (\Ge^{-1} X)| = \Gb_j |\Ge^{-1} X-S_j|^{\Gb_j-1} \simeq \Ge^{-\Gb_j+1} |X-V_j|^{\Gb_j-1} , \quad j=1,2,
\eeq
we infer from \eqnref{B1B2est} that
\beq
\big| (\nabla \Bcal_1 - b \nabla \Bcal_2) ( \Ge^{-1} X ) \big| \simeq \sum_{j=1}^2 |(\nabla \Bcal_j ) ( \Ge^{-1} X )| \simeq \sum_{j=1}^2 \Ge^{-\Gb_j+1} |X-V_j|^{\Gb_j-1}.
\eeq
As one can see from \eqnref{error}, we treat the term $\Ncal(\Ge^{-1}X)$ as an error term. Since $\Ncal(\Ge^{-1}X)$ is bounded, the decomposition \eqnref{naurep} is meaningful only if $|(\nabla \Bcal_1 - b \nabla \Bcal_2) ( \Ge^{-1} X )|$ is bounded below, or $|X-V_j| \le c_0 \Ge$ for some $c_0$. In such a case, we have
\begin{align*}
\big| (\nabla \Bcal_1 - b \nabla \Bcal_2) ( \Ge^{-1} X ) + \Ncal(\Ge^{-1}X) \big|
& \ge \big| (\nabla \Bcal_1 - b \nabla \Bcal_2) ( \Ge^{-1} X ) \big| - \big| \Ncal(\Ge^{-1}X) \big| \\
& \gtrsim \sum_{j=1}^2 \Ge^{-\Gb_j+1} |X-V_j|^{\Gb_j-1}.
\end{align*}
Furthermore, if $|X-V_k| \le c_0 \Ge$, then
$$
\sum_{j=1}^2 \Ge^{-\Gb_j+1} |X-V_j|^{\Gb_j-1} \simeq \Ge^{-\Gb_k+1} |X-V_k|^{\Gb_k-1}.
$$
So, we obtain the following theorem.
\begin{thm}\label{1st_cor}
\begin{itemize}
\item[(i)] The following holds for $X \in B_{\Gd_1} \setminus \overline{\GO_1 \cup \GO_2}$
\beq\label{naurep3}
|\nabla u(X)| \lesssim \frac{\norm{ h}_{L^{\infty} (\GO)}}{\Ge|\log \Ge|} \left[ \sum_{j=1}^2 \Ge^{1-\Gb_j} |X-V_j|^{\Gb_j-1} + 1 \right].
\eeq
\item[(ii)] Suppose that $|a_\Ge| \ge a_0$ for some positive constant $a_0$. There is a constant $c_0$ independent of $\Ge$ and $X$ such that
\beq\label{naurep4}
|\nabla u(X)| \gtrsim \frac{\Ge^{-\Gb_j}}{|\log \Ge|} |X-V_j|^{\Gb_j-1} \quad\mbox{for all } X \in B_{c_0\Ge}(V_j), \quad j=1,2.
\eeq
\end{itemize}
\end{thm}

Inequalities \eqnref{naurep3} and \eqnref{naurep4} show that the corner singularity  $|X-V_j|^{\Gb_j-1}$ is amplified by the factor of $\Ge^{-\Gb_j}/|\log \Ge|$, which is due to the interaction between two perfectly conducting inclusions. It is worth emphasizing that the blow-up magnitude $\frac{\Ge^{-\Gb_j}}{|\log \Ge|} |X-V_j|^{\Gb_j-1}$ is much larger than that for the case when inclusions have smooth boundaries, namely, $\Ge^{-1/2}$, as mentioned in Introduction.

\section{A lower estimate of $a_\Ge$}\label{sec:aGe}

The estimates from below in Theorems \ref{cortwo} and \ref{1st_cor} are obtained under the assumption that $|a_\Ge| \ge a_0$ for some positive constant $a_0$ independent of $\Ge$. We now show that this assumption can be fulfilled for some background potential $h$. It is worth emphasizing that there are other cases of $h$ where $a_\Ge=0$. For example, if inclusions $\GO_1$ and $\GO_2$ are symmetric with respect to the $y$-axis and $h(x,y)=y$, then one can easily see that $u(-x,y)$ is a solution to \eqnref{main} if $u(x,y)$ is. It implies that $u|_{\p\GO_1}=u|_{\p\GO_2}$, and hence $a_\Ge=0$.

We have the following theorem (without assuming that $\GO_1$ and $\GO_2$ are symmetric with respect to the $y$-axis).

\begin{thm}\label{2nd_thm}
If $h(x,y) = x$, then
\beq
a_{\Ge} \simeq 1.
\eeq
\end{thm}

\pf
Recall that
$$
a_{\Ge} = \frac{\Gg |\log\Ge|}{\pi (q |_{\p\GO_2} - q|_{\p\GO_1})} (q |_{\p \GO_2} - q|_{\p \GO_1})(u |_{\p \GO_2} - u|_{\p \GO_1})  .
$$
It is already proved in \eqnref{aGe} that $a_\Ge$ is bounded above regardless of $\Ge$ .

By Lemma \ref{lemma_C}, there exists a constant $C_1$ such that
\beq\label{aGelower}
\frac{\Gg |\log\Ge|}{\pi (q |_{\p\GO_2} - q|_{\p\GO_1})} > C_1 >0.
\eeq
On the other hand, it follows from \eqnref{q_and_Q} and \eqref{the_presentation_potential_diff_u} that
\begin{align*}
(q |_{\p \GO_2} - q|_{\p \GO_1})(u |_{\p \GO_2} - u|_{\p \GO_1}) &= (q |_{\p \GO_2} - q|_{\p \GO_1}) \int_{\p \GO_1 \cup \p \GO_2} x \p_{\nu} q ds \\
&= \int_{\p \GO_1 \cup \p \GO_2} x \p_{\nu} \Gn ds.
\end{align*}
Since $(-1)^{i} \p_{\nu} q > 0 $ on $\p \GO_i$ for $i=1,2$, we have
$$
(q |_{\p \GO_2} - q|_{\p \GO_1})(u |_{\p \GO_2} - u|_{\p \GO_1}) \ge \int_{\p\GO_2} x \p_{\nu} \Gn ds.
$$

To estimate $\int_{\p \GO_2} x \p_{\nu} \eta ds$, we choose positive constants $r_0$ and $\Gd_0$ regardless of small $\Ge>0$ so that
$$
B_{r_0} (V_1-(\Gd_0,0)) \subset \GO_1 .
$$
Let $B_0:= B_{r_0} (V_1-(\Gd_0,0))$ for ease of notation. Let $\Gz$ be the solution to
$$
\GD \Gz = 0 \quad\mbox{in }  \Rbb^2 \setminus \overline{B_{0}\cup \GO_2}
$$
with conditions
$$
\begin{cases}
\ds \Gz = 0  \quad\mbox{on } \p B_{0} ,  \\
\ds \Gz = 1  \quad\mbox{on } \p \GO_{2}  ,\\
\ds \int_{\Rbb^2 \setminus \overline {B_{0}\cup \GO_2}} |\nabla \Gz|^2 dxdy < \infty .
\end{cases}
$$
We emphasize that $\Gz$ is independent of $\Ge$. Since
$$
1 \geq \Gz(X) \geq \Gn(X) \geq 0 \quad\mbox{for } X \in \Rbb^2 \setminus \overline{\GO_1 \cup \GO_2}
$$
and
$$
\Gz= \Gn = 1 \quad\mbox{on }\p \GO_2,
$$
we infer from Hopf's lemma that
$$
0 < \p_{\nu} \Gz \leq \p_{\nu} \Gn \quad\mbox{on }\p \GO_2.
$$
Hence, there exists a constant $C_2>0$ regardless of $\Ge$ such that
\begin{align*}
\int_{\p \GO_2} x \p_{\nu} \Gn ds \geq \int_{\p \GO_2} x  \p_{\nu} \Gz ds \geq C_2.
\end{align*}
This together with \eqnref{aGelower} shows that $a_\Ge$ is bounded below regardless of $\Ge$, and the proof is complete.
\qed

\section*{Conclusion}

In this paper we deal with the field enhancement due to presence of a bow-tie structure consisting of two perfectly conducting inclusions with corners. We introduce auxiliary functions to capture singular behavior of the field as the distance $\Ge$ between two inclusions tends to $0$. As consequences we obtain optimal bounds on the size of the field in the region close to vertices and the region relatively away from (but still close to) vertices. The estimates show that the field is enhanced beyond the corner singularities due to the interaction between inclusions. There may be a gap between two regions of estimates. In the gap the lower bound on $|\nabla u|$ is not obtained, and it is interesting to fill the gap.  

\section*{Acknowledgement}
We thank S. Yu for pointing out existence of the reference \cite{PBFLN}.


\begin{thebibliography}{99}

\bibitem{ACKLY-ARMA-13} H. Ammari, G. Ciraolo, H. Kang, H. Lee and K. Yun, Spectral analysis of the Neumann-Poincar\'e operator and characterization of the stress concentration in anti-plane elasticity, Arch. Ration. Mech. Anal. 208 (2013), 275--304.

\bibitem{AKL-MA-05} H. Ammari, H. Kang and M. Lim, Gradient estimates for solutions to the conductivity problem, Math. Ann. 332(2) (2005), 277--286.

\bibitem{bab} I. Babu\u{s}ka, B. Andersson, P. Smith and K. Levin, Damage analysis of fiber composites. I. Statistical
analysis on fiber scale, Comput. Methods Appl. Mech. Engrg. 172 (1999), 27--77.

\bibitem{BLL-ARMA-15} J. Bao, H. Li and Y. Li, Gradient estimates for solutions of the Lam\'e
system with partially infinite coefficients, Arch. Rational Mech. Anal. 215 (2015), 307--351.

\bibitem{BLY-ARMA-09}  E.S. Bao, Y. Li and B. Yin,
Gradient estimates for the perfect conductivity problem, Arch. Rat. Mech. Anal. 193 (2009), 195-226.

\bibitem{BLY-CPDE-10} E.S. Bao, Y. Li and B. Yin,
Gradient estimates for the perfect and insulated conductivity problems with multiple inclusions, Commun. Part. Diff. Eq. 35 (2010), 1982--2006.

\bibitem{DL-arXiv} H. Dong and H. Li, Optimal estimates for the conductivity problem by Green's function method, arXiv:1606.02793.

\bibitem{GN-MMS-12} Y. Gorb and A. Novikov, Blow-up of solutions to a $p$-Laplace equation, SIAM Multi. Model. Simul. 10 (2012), 727--743.

\bibitem{Grisvard-book} P. Grisvard, {\sl Boundary value problems in non-smooth domains}, Pitman,
London, 1985.

\bibitem{KLY-MA-15} H. Kang, H. Lee and K. Yun, Optimal estimates and asymptotics for the stress concentration between closely located stiff inclusions, Math. Annalen 363 (2015), 1281--1306.

\bibitem{KLY-JMPA-13} H. Kang, M. Lim and K. Yun, Asymptotics and computation of the solution to the conductivity
equation in the presence of adjacent inclusions with extreme conductivities, {J. Math. Pure. Appl.} 99 (2013), 234--249.

\bibitem{KLY-SIAP-14} H. Kang, M. Lim and K. Yun, Characterization of the electric field concentration between two adjacent spherical perfect conductors, SIAM J. Appl. Math. 74 (2014), 125--146.

\bibitem{KY-arXiv} H. Kang and S. Yu, Qualitative characterization of stress concentration in presence of adjacent hard inclusions in two-dimensional linear elasticity, preprint.

\bibitem{Kondra-TMMS-67} V. A. Kondratiev, Boundary-value problems for elliptic equations in
domains with conical or angular points. Trans. Moscow Math. Soc. 16 (1967), 227--313.

\bibitem{KMR-book} V.A. Kozlov, V.G. Maz'ya and J. Rossmann, {\sl Elliptic boundary
value problems in domains with point singularities}, Amer. Math. Soc., Mathematical Surveys and Monographs, vol 52, Providence, RI, 1997.

\bibitem{PBFLN} V. Pacheco-Pen\~{a}, M. Beruete, A.I. Fern\'{a}ndez-Dom\'{i}nquez, Y. Luo and M. Navarro-C\'{i}a, Description of bow-tie nanoantennas excited by localized emitters using conformal transformation, ACS Photonics 2016, 3, 1223−-1232.

\bibitem{Y} {K. Yun,} Estimates for electric fields blown up between closely adjacent conductors with arbitrary shape, {SIAM J. Appl. Math.} { 67}  (2007), 714--730.

\bibitem{Yun-JDE-16} K. Yun, An optimal estimate for electric fields on the shortest line segment between two spherical insulators in three dimensions, J. Differ. Equations 261 (2016), 148-188.

\end{thebibliography}
\end{document}